\documentclass{elsart}

\usepackage{amsmath}
\usepackage{amssymb}
\usepackage{amsfonts}
\usepackage{amscd}
\usepackage{graphicx} 
\usepackage{epsfig}

\newtheorem{THM}{Theorem}[section]
\newtheorem{PROP}[THM]{Proposition}

\newtheorem{LEM}[THM]{Lemma}
\newtheorem{CLAIM}[THM]{Claim}

\begin{document}
\begin{frontmatter}


\title{The Hexatangle}


\author{Lorena Armas-Sanabria\thanksref{label2}},
\author{Mario Eudave-Mu\~noz\corauthref{cor1}\thanksref{label2}\thanksref{label3}}
\ead{larmas@matem.unam.mx}
\ead{mario@matem.unam.mx}
\thanks[label2]{Instituto de Matem\'aticas, Universidad Nacional Aut\'onoma de M\'exico, Ciudad Universitaria, 04510 M\'exico D.F., M\'exico}
\thanks[label3]{CIMAT, Callej\'on Jalisco s/n, Guanajuato, Gto. M\'exico}

\corauth[cor1]{To whom the correspondence should be addressed}

\address{}

{\small\it Dedicated to Michel Domergue on the occasion of his sixtieth birthday and to the memory of Yves Mathieu}

\begin{abstract}
We are interested in knowing what type of manifolds are obtained by doing
Dehn surgery on closed pure $3$-braids in $S^3$. In particular, we want to
determine when we get $S^3$ by surgery on such a link. 
We consider links which are small closed pure 3-braids; these are the
closure of 3-braids of the form 
$({\sigma_1}^{2e_1})({\sigma_2}^{2f_1})(\sigma_2\sigma_1\sigma_2)^{2e}$, where $\sigma_1$, 
$\sigma_2$ are
the generators of the 3-braid group and $e_1$, $f_1$, $e$ are integers.
We study Dehn surgeries on these links, and determine exactly which ones
admit an integral surgery producing the 3-sphere. This is
equivalent to determining the surgeries of some type on a certain six component link $\mathcal{L}$ that produce $S^3$. The link $\mathcal {L}$
is strongly invertible and its exterior double branch covers a certain
configuration of arcs and spheres, which we call the Hexatangle. Our problem is
equivalent to determine which fillings of the spheres by integral tangles
produce the trivial knot, which is what we explicitly solve. 
This hexatangle is a generalization of the
Pentangle, which is studied in \cite{GordonL3}.
\end{abstract}

\end{frontmatter}

\section{Introduction}
\label{sec:intro}

We are interested in knowing what type of
manifolds are obtained by doing Dehn
surgery on closed pure 3-braids in $S^3$. In
particular, when is possible to obtain the 3-sphere by 
Dehn surgery on a closed pure 3-braid.

By the fundamental theorem of surgery proved by
Lickorish and Wallace \cite{Lickorish1}, \cite{Lickorish2}, \cite{Wallace},
we know that any closed, connected and oriented
3-manifold can be obtained by integral Dehn surgery
on a closed pure $n$-braid. It is known that surgery on a
closed pure 1-braid produces lens spaces, for such
a braid is the trivial knot;
some surgeries
on closed pure 2-braids produce connected sums of
lens spaces, but in general they produce
Seifert fibered spaces, for a closed pure $2$-braid is
a torus link. 
So, it is a natural
question to ask  what kind of 3-manifolds are obtained by
surgery on closed pure 3-braids.

By \cite{FadellN} we have that the group of pure 3-braids can be seen
as the direct product of two free groups $Z \times F_2$. 
So the group of pure $3$-braids can be expressed as 
$P_3 = \{  \beta\in B_3  \ \vert \ \beta = \Delta^{2e} 
\prod \sigma _1^{2e_i } \sigma _2^{2f_i} \} $, where 
$\Delta = \sigma_2\sigma_1\sigma_2$, and $e,\ e_i,\ f_i$ are integers. 
Denote by $\hat \beta$ the closure of the braid $\beta$.

In this work we consider the closed 3-braids of the form 
$\hat \beta = \widehat {\sigma_1^{2e_1}\sigma_2^{2f_1}({\sigma_2\sigma_1\sigma_2})^{2e}}$,
shown in Figure 1, where the boxes indicate the number of full twist given to the braid.
We call these links {\it small} closed pure 3-braids.
We determine when an integral surgery in such a link produces the 3-sphere.

In a previous work \cite{ArmasE}, we considered closed pure
3-braids $\hat \beta$ of the form
$\hat \beta = \widehat { \prod _{i=1 }^n \sigma_1^{2e_i}\sigma_2^{2f_i }}$ where
$\vert e_i \vert \geq 1$, $\vert f_i \vert \geq 1$, $n\geq 2$, and showed that 
in many cases we obtain a Haken or a laminar manifold by
surgery on such links.
The first author \cite{Armas} has shown an example
of a hyperbolic small closed pure $3$-braid which do have a nontrivial surgery
producing $S^3$, which is recovered in the present paper.

\begin{figure}[htbp]
\begin{center}
\includegraphics[width=0.3\linewidth]{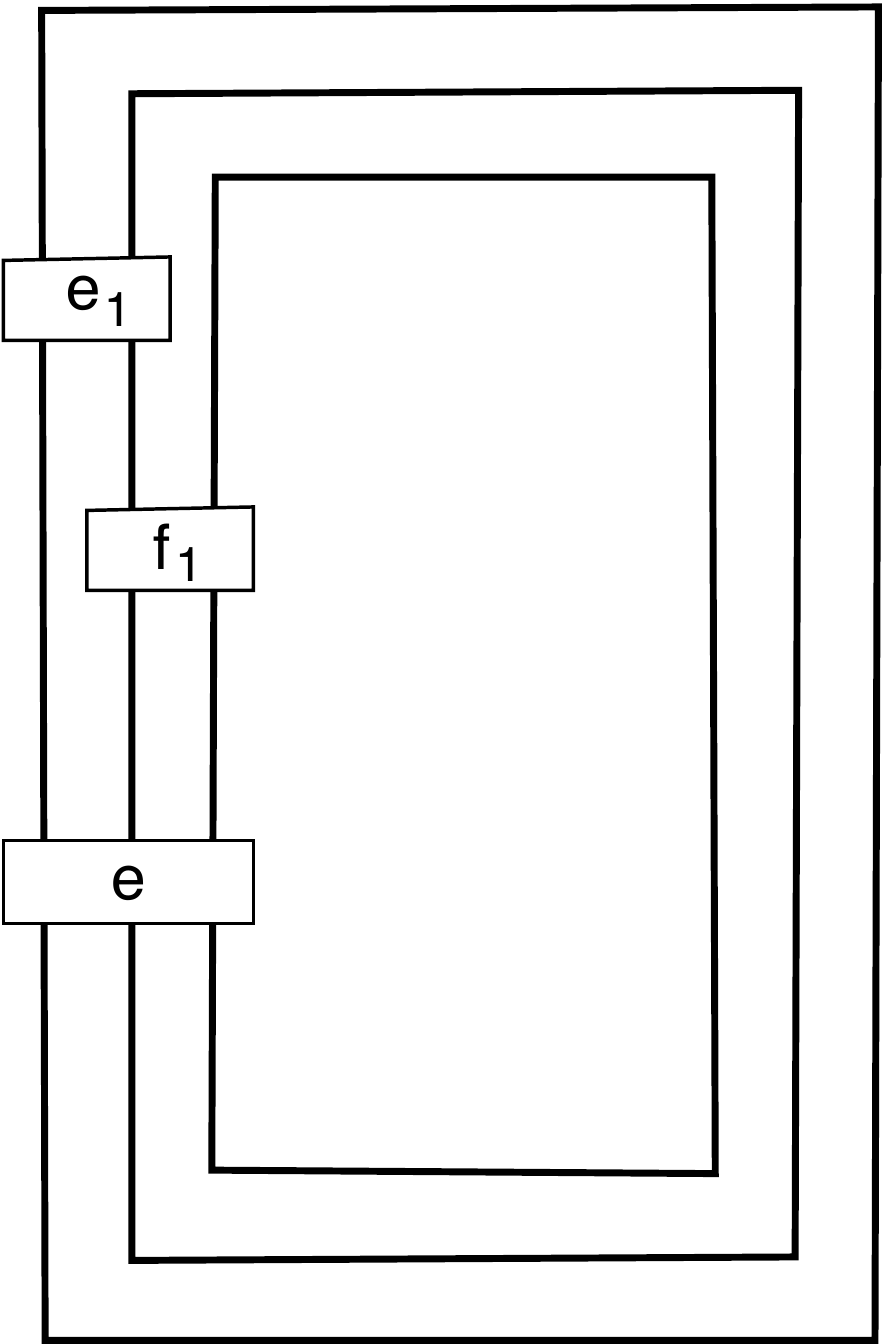}
\caption{}
\label{Figure1}
\end{center}
\end{figure}

Note that the link $\hat \beta$ which is the closure of the braid
${\sigma_1^{2e_1}\sigma_2^{2f_1}({\sigma_2\sigma_1\sigma_2})^{2e}}$, can be obtained by
$(1/e_1,1/f_1,1/e)$-Dehn surgery on the link $\mathcal {L}$ shown in Figure 2.
It is known that this link is hyperbolic, and in fact arithmetic \cite{Baker}.
So our problem is equivalent to determine when surgery on this link produces the
3-sphere, though we consider integral surgery on 3 components of the link 
and integral reciprocal in the other 3 components. We indicate surgeries on this link
by $\mathcal {L}(1/e_1,1/f_1,1/e,m,n,p)$, as indicated in Figure 2, which implicitly is
giving an order to the components of the link.

\begin{figure}[htbp]
\begin{center}
\includegraphics[width=0.6\linewidth]{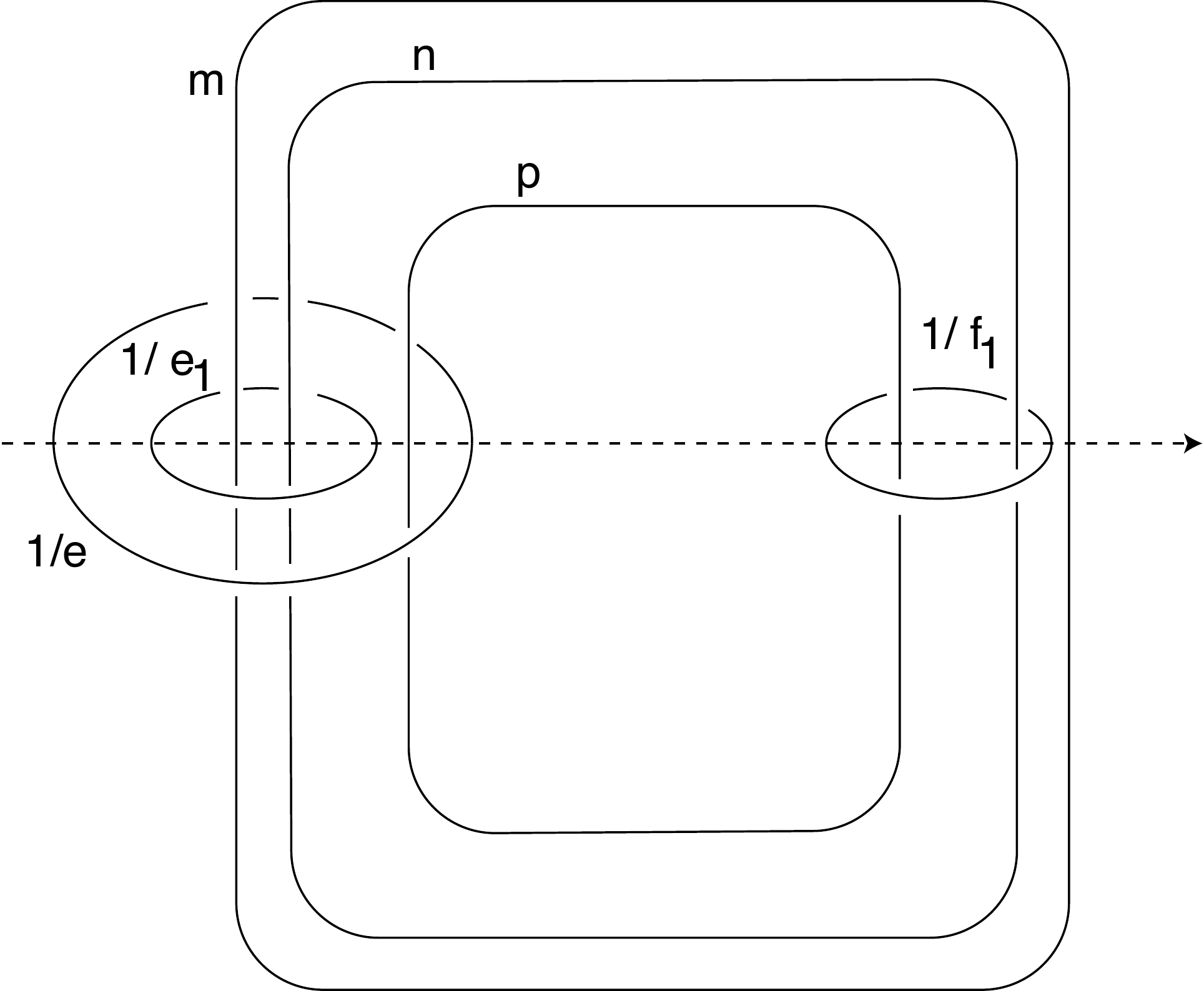}
\caption{}
\label{Figure2}
\end{center}
\end{figure}

Note that the link $\mathcal {L}$ is strongly invertible, an involution axis is shown in 
Figure 2. The quotient of the exterior of $\mathcal {L}$ under this involution will be a 
punctured $S^3$, together with arcs, formed by the image of the involution axis.

More precisely, following \cite{GordonL4}, a {\it tangle} will be a pair $(B,A)$ where 
$B$ is $S^3$ with the interiors of a finite number ($\geq 1$) of disjoint 3-balls 
removed, and $A$ is a disjoint union of properly embedded arcs in $B$ such that $A$ 
meets each component of $\partial B$ in exactly four points. Two tangles 
$(B_1,A_1)$ and 
$(B_2,A_2)$ are homeomorphic 
if there is a homeomorphism of pairs $h:(B_1,A_1)\rightarrow (B_2,A_2)$. 

A {\it marking} of a tangle $(B,A)$ is an identification of each pair 
$(S,S\cap A)$, where $S$ is a component of $\partial B$, with 
$(S^2,Q=\{ NE,NW,SW,SE\})$. A {\it marked} tangle is a tangle together with a marking. 
We say that a homeomorphism preserves a marking if the axis
$NW-NE$ is mapped to one of the axes $NW-NE$, $NE-NW$, $SW-SE$ or $SE-SW$, and the 
other axes are mapped accordingly. 
Two marked tangles are equivalent if they are homeomorphic by an orientation preserving 
homeomorphism that preserves the markings.

A {\it rational} tangle is a marked tangle that is homeomorphic to the trivial tangle 
in the 3-ball, $(D^2, 2\ points)\times I$. As marked tangles, rational tangles are 
parameterized by $Q\cup \{1/0\}$. We denote the rational tangle corresponding to
$p/q\in Q\cup \{1/0\}$ by $R(p/q)$, and adopt the conventions of \cite{Eudave3}. Given a marked tangle, there is a well defined way of filling its boundary components with rational tangles. 
 
So, the quotient of $\mathcal {L}$ under the involution is a tangle (see Figure 3), where 
its boundary 
components come from the tori boundary components of the exterior of $\mathcal {L}$, and 
the arcs are the image of the involution axis. This tangle could have a natural 
marking, if we choose it as given by the image of a framing on the components 
of $\mathcal {L}$, 
as shown in Figure 3. Instead we choose a marking as in Figure 4. This is indicated 
in Figure 4 by a rectangular box, where the short sides of the rectangle represent the 
axis $NW-SW$ and $NE-SE$, and the long sides represent the axis $NW-NE$ and $SW-SE$. In all of our pictures the shape of the rectangle will be always clear.
We call this marked tangle the {\it Hexatangle}, and denote it by $\mathcal {H}$,
or $\mathcal {H}(*,*,*,*,*,*)$. The capital letters $A,B,C,D,E,F$ denote boundary components in the
hexatangle, and $\alpha,\beta,\gamma,\delta,\epsilon,\eta$ denote fillings of the hexatangle 
with rational tangles, so for example, $\mathcal {H}(\alpha,\beta,*,*,*,\eta)$ denote the 
tangle obtained by filing the components $A$, $B$ and $F$ with the rational tangles 
$\alpha$, $\beta$ and $\eta$ respectively. 
We call the sphere boundary components of $\mathcal {H}$, filled or unfilled, simply boxes. We say that
two  boxes are adjacent if there is an arc of $\mathcal {H}$ connecting them, and opposite otherwise.
So each box is opposite to just one box and adjacent to 4 boxes.
We consider $\alpha,\beta,\gamma,\delta,\epsilon,\eta$
as rational parameters, so that when $\alpha=-1$, we mean that we are filling the
corresponding box with the integral tangle $R(-1)$. Note that when we fill the boxes with integral tangles, we are just replacing each box with a sequence of horizontal crossings.

We remark that the hexatangle is the same as Conway's basic polyhedra $6^*$ and $6^{**}$,
but with different marking \cite{Conway}, \cite{Jablan}.

\begin{figure}[htbp]
\begin{center}
\includegraphics[width=0.9\linewidth]{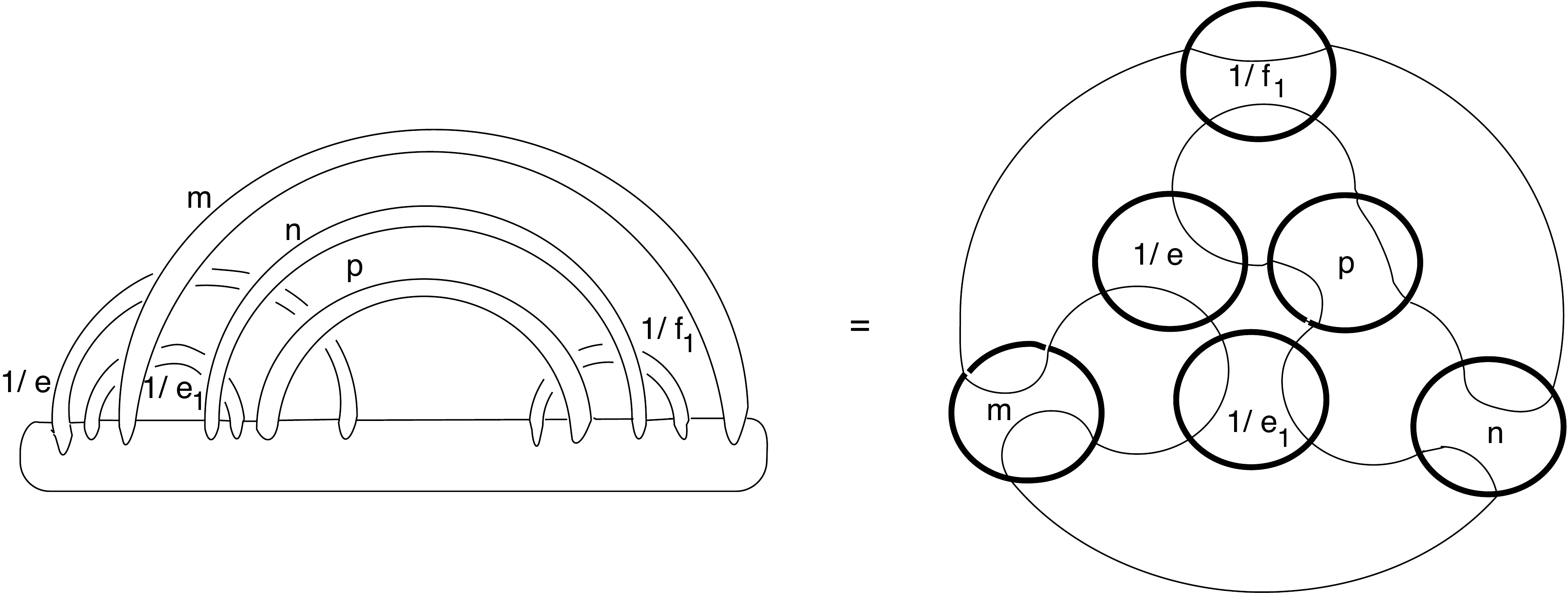}
\caption{}
\label{Figure3}
\end{center}
\end{figure}

\begin{figure}[htbp]
\begin{center}
\includegraphics[width=0.8\linewidth]{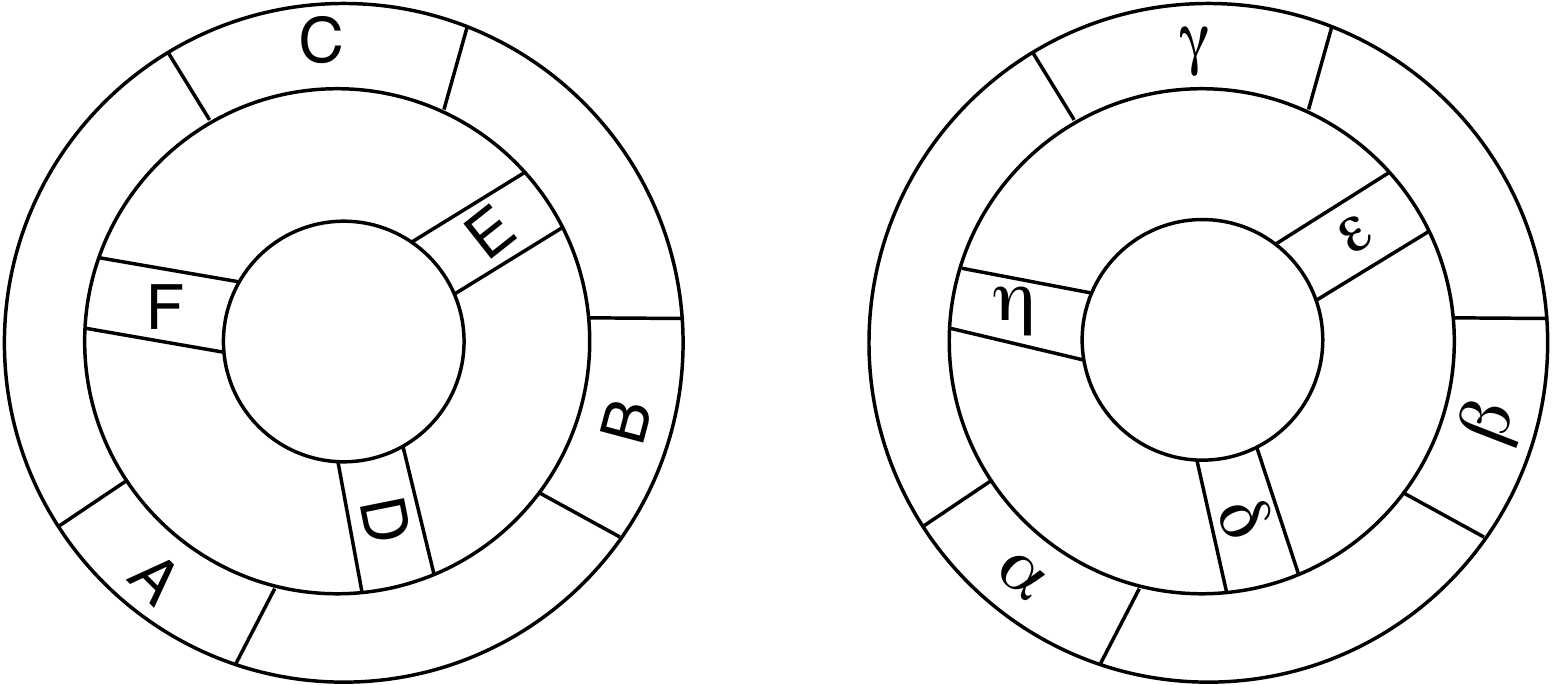}
\caption{}
\label{Figure4}
\end{center}
\end{figure}

Note that by filling one of the components $A$, $B$, $E$, with a rational tangle 
$\mathcal {R}(p/q)$ will correspond in the double branched cover, to do $(-p/q)$-Dehn 
surgery on the corresponding component, while filling with $\mathcal {R}(p/q)$ in one of the 
components $C$, $D$, $F$, will correspond in the double branched cover, 
to do $q/p$-Dehn surgery in the corresponding component (see \cite{Montesinos}), this because
of our rational tangles convention (see \cite{Eudave3}). So we can consider integral fillings
in all boundary components of the hexatangle and forget the correspondence 
with the components of $\mathcal {L}$.

Remember that the 3-sphere double branch covers only the trivial knot, by the
solution of the Smith conjecture. 
So, our original problem about surgery on small closed pure 3-braids translate to the
following:

{\it When is it possible to get a trivial knot by filling the hexatangle with 
integral tangles ?}

The same question could be asked for any fillings, i.e., when the trivial knot is 
obtained by rational fillings of the hexatangle? Determining this is equivalent to 
determining all Dehn surgeries on the link $\mathcal {L}$ that produce the 3-sphere.
We plan to study this problem in a subsequent paper.

This problem is interesting by itself but also for several other reasons.
Lots of hyperbolic manifolds of small volume are obtained by doing surgery on some 
components of $\mathcal {L}$, for example, by doing $(-1)$-surgery on any component of 
$\mathcal{L}$, we get a link whose exterior is isometric to the exterior of the minimally. Also, the Pentangle which is studied in \cite{GordonL3}, is obtained by
putting $\beta=1$ in the Hexatangle.  The so called ``magic manifold'',
which is the exterior of the 3-chain link studied in \cite{MartelliP}, is also
obtained by Dehn surgery on $\mathcal {L}$. In fact, the 3-chain link is the closure the pure 3-braid
${\sigma_1^{4}\sigma_2^{4}({\sigma_2\sigma_1\sigma_2})^{-2}}$. In
\cite{MartelliP} all the  exceptional fillings of the 3-chain link are determined; these results can be verified 
by looking at the corresponding fillings of the hexatangle. It would be also interesting 
to determine all exceptional fillings of the link $\mathcal {L}$.

D. Futer and J.S. Purcell \cite{FPurcell} have shown that if a link $K$ has a prime, twist-reduced diagram $D(K)$, with at least two twist regions and each twist region containing at least 6 crossings, then $K$ is hyperbolic. This implies that by filling the hexatangle with integral tangles, each in absolute value greater or equal to 6, then we get hyperbolic links, in particular the trivial knot is not obtained.
Here we give a sharp result for the hexatangle, showing exactly when we get the trivial knot. 
It would also be interesting to determine when a non-hyperbolic link is obtained from the hexatangle.

\vfill\eject

Another reason why it is interesting to determine when we get the trivial knot by 
filling the hexatangle, is that if a certain filling produce the trivial knot, then 
by filling all the components except one, we get a 2-string tangle whose double 
branched cover is the exterior of a knot in $S^3$, or in other words, by doing surgery 
on the corresponding five components of $\mathcal {L}$, we get the exterior of a knot in $S^3$. 
By experimentation, we can see that many of those knots are hyperbolic and have 
non-hyperbolic surgeries, in fact, Seifert fibered space surgeries. Many of the examples that we know of hyperbolic strongly invertible knots with a Seifert fibered surgery, come from surgery on $\mathcal{L}$, ref. \cite{Eudave3}, \cite{Tera1}, \cite{Baker2}. So, by solving the problem about the hexatangle we could get an interesting list of hyperbolic strongly invertible knots with a Seifert fibered space surgery.
However, we cannot expect to find all hyperbolic strongly invertible knots with a Seifert fibered space in this way, for the volume of a knot with a lens space surgery can be arbitrarily large \cite{Baker1}, while the volume of any hyperbolic knot obtained by surgery on $\mathcal{L}$ is bounded. We remark that there are hyperbolic non-strongly invertible knots with Seifert 
fibered surgeries \cite{MMM},\cite{DMM}, \cite{Tera2}. Also, many examples of hyperbolic manifolds with
exceptional fillings constructed via tangles, are special cases of the hexatangle
(ref. \cite{EudaveWu}).

The hexatangle has many symmetries. Note that the hexatangle can be embedded in a tetrahedron,
so that each box is in correspondence with an edge of the tetrahedron, as shown in Figure 5. So any symmetry of the tetrahedron will give a symmetry of the hexatangle preserving framings. We give a list of fillings on the hexatangle that produces that produces the trivial knot up to
symmetries, where by this we mean that the list is complete up to  the symmetries given by the tetrahedron and mirror images.
Note that given any two boxes, there is a symmetry that
takes one to the other. Also, given two pairs of adjacent
(opposite) boxes there is a symmetry that takes one pair to the other.

\begin{figure}[htbp]
\begin{center}
\includegraphics[width=0.45\linewidth]{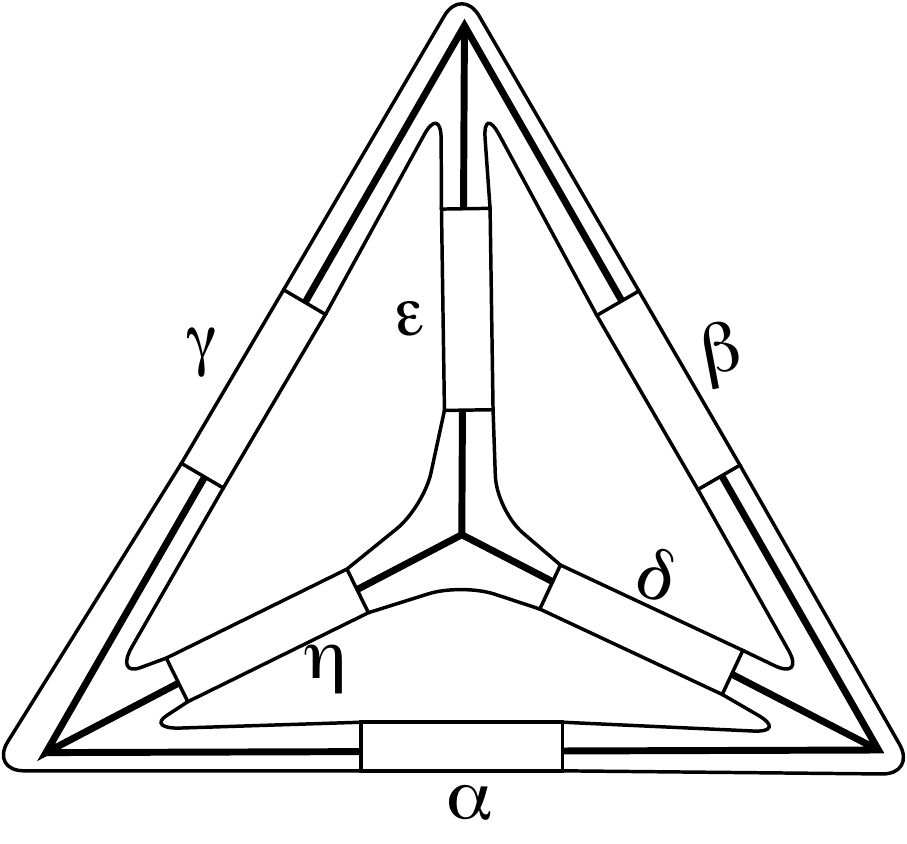}
\caption{}
\label{Figure5}
\end{center}
\end{figure}

\vfill\eject

 Our results are the following:

\begin{THM}\label{THM:main}
Suppose an integral  filling of the hexatangle produces the trivial knot,
then the parameters are exactly as shown in Tables 1, 2 and 3, up to symmetries.
\end{THM}

This will follows from the following results.

\begin{THM}\label{THM:C-0}
Suppose that one of the parameters 
$\alpha$, $\beta$, $\gamma$, $\delta$, $\epsilon$, $\eta$,
say $\eta$ is the tangle $0$. Then 
$\mathcal {H} (\alpha,\beta,\gamma,\delta, \epsilon,0)$ is the trivial
knot if and only the parameters are as in Tables 1 and 2, up to symmetries.
\end{THM}

\begin{THM}\label{THM:quasi-main} 
Suppose that all of 
$\alpha$, $\beta$, $\gamma$, $\delta$,
$\epsilon$ and $\eta$ are different from $0$. If 
$\mathcal {H} (\alpha,\beta,\gamma,\delta, \epsilon,\eta)$ is the trivial
knot, then there is a pair of adjacent boxes, say $\delta$ and $\eta$,
so that $\delta=-1$ and $\eta=1$.
\end{THM}

\begin{THM}\label{THM:C-1}
 Suppose that $\alpha$, $\beta$, $\gamma$ and
$\epsilon$ are not $0$, $\delta=-1$ and $\eta=1$.
Then $\mathcal {H} (\alpha,\beta,\gamma,-1, \epsilon,1)$ is the trivial
knot if and only the parameters are as in Table 3, up to symmetries.
\end{THM}

Theorems \ref{THM:C-0} and \ref{THM:C-1} are just a rational tangles computation; this is carried out
in Sections 2 and 3, respectively. Section 4 is dedicated to a proof of Theorem \ref{THM:quasi-main};
this is the main part of the paper. First, we apply some deep results on Dehn surgery 
on knots to reduce the theorem to six cases, in which there are two small boxes
(Lemma \ref{LE:6C}). Then an analysis is made of each of the cases, to conclude that
the trivial knot cannot be obtained, except in one of the cases.
In Section 5 we discuss about surgery on closed pure 3-braid producing $S^3$,
and show that there are infinitely many hyperbolic small closed pure 3-braid with 
a nontrivial surgery producing the 3-sphere.

\section{When a parameter is 0}
\label{sec:paremeter0}

In this and next section we do some rational tangles computations and rely on
known facts about rational tangles and knots. We follow the conventions of
\cite{Eudave3}.
We denote by $R(p/q)$ the rational tangle determined by $p/q$, and by $K(p/q)$
the rational knot or 2-bridge knot, which is the numerator of the rational tangle
$R(p/q)$. As usual, the numerator of a rational tangle is obtained by closing it
with two arcs, one arc joining the points NW-NE, and the other the points
SW-SE, without introducing new crossing. A rational tangle
$R(p/q)$ can be given by a sequence of integers
$[a_1,a_2,\dots,a_n]$ whose continued fraction gives $p/q$, 
i.e. $p/q = a_n + {1\over{a_{n-1} + {1\over{a_{n-2}+ ...}}}}$; in this case
$R[a_1,a_2,\dots,a_n]$ denotes the tangle $R(p/q)$ and $K[a_1,a_2,\dots,a_n]$
denotes the numerator of such tangle. 
If $K(p/q)$ is the trivial knot, and $p$, $q$ are relative primes, then $p=\pm 1$.
Also note that if $p/q$ is obtained from a continued fraction, then in fact
$p$ and $q$ are relative primes.

A Montesinos tangle $T(p_1/q_1,\dots,p_n/q_n)$ is a tangle formed by a horizontal strand of
rational tangles, and a Montesinos link $M(p_1/q_1,\dots,p_n/q_n)$ is the numerator 
of a Montesinos tangle. The double cover of
$B^3$ branched along a Montesinos tangle $T(p_1/q_1,\dots,p_n/q_n)$ is a Seifert fibered
space $D(p_1/q_1,\dots,p_n/q_n)$ over the disk with at most $n$-cone points of index
$q_1,\dots,q_n$. The double cover of $S^3$ branched along
$M(p_1/q_1,\dots,p_n/q_n)$ is a Seifert fibered space over the sphere $S^2$ with 
at most $n$-cone points
of index $q_1,\dots,q_n$. So, if $M(p_1/q_1,p_2/q_2,p_3/q_3)$ is a trivial knot,
one of the tangles $R(p_i/q_i)$ is an integral tangle, so that it can be inserted
into one of the other rational tangles, getting a Montesinos knot formed by two tangles,
that is, a 2-bridge knot. Note also that if $M(p_1/q_1,p_2/q_2,p_3/q_3)$ is a composite 
link, then one of the tangles $R(p_i/q_i)$ is the rational tangle $R(1/0)$.
Finally note that if the trivial knot is presented as a sum of two 2-strings tangles, then
at least one of the tangles must be a trivial tangle. In what follows we use the words 
knot and link interchangeably, to mean a collection  of circles, except when referring to 
the trivial knot, which always will consist of a single component.

In this section we prove
the following,

{\bf Theorem 1.2}\ {\it Suppose that one of the parameters 
$\alpha$, $\beta$ $\gamma$, $\delta$, $\epsilon$, $\eta$,
say $\eta$ is the tangle $0$. Then 
$\mathcal {H} (\alpha,\beta,\gamma,\delta, \epsilon,0)$ is the trivial
knot if and only if the parameters are as in Tables 1 and 2, up to symmetries.}

\textbf{Proof} The proof is a rational tangles calculation. Note that if 
4 or more of the parameters are $0$, then the link obtained has more than one
component. So suppose at most 3 of the parameters are $0$.

If 3 of the parameters are $0$, then we have two cases up to symmetry:
A) $\delta=0$, $\eta=0$, $\epsilon=0$, and B) $\beta=0$, $\delta=0$, $\eta=0$.
All other possible cases would produce a link of several components.

Case A. $\delta=0$, $\eta=0$, $\epsilon=0$

Here the knot looks like a connected sum of 3 knots, so it can be the trivial knot 
if an only if $\alpha=\pm 1$,
$\beta=\pm 1$ and $\gamma=\pm 1$. This makes line 1 of Table 1.

Case B. $\beta=0$, $\delta=0$, $\eta=0$

Here the knot can be the trivial knot if an only if $\alpha=\pm 1$, $\gamma=\pm 1$
and $\epsilon=\pm 1$. This makes line 2 of Table 1.

\vfill\eject

Suppose now that just two of the parameters are $0$. Here we have two cases, 
up two symmetries, depending if the given boxes are adjacent or opposite:
C) $\delta=0$, $\eta=0$ and D) $\beta=0$, $\eta = 0$.

Case C. $\delta=0$, $\eta=0$

In this case the knot looks like a composite knot, so to be trivial both 
components must be trivial. One
of them is trivial if and only if $\alpha=\pm 1$. The other one is the Montesinos knot
$M(-1/\gamma,-1/\epsilon,-1/\beta)$, so to be trivial one of $\gamma$, 
$\epsilon$ or $\beta$ must be  $\pm 1$. 

Case C.1. $\beta =1$

In this case, the knot is isotopic to the rational knot
$K[\gamma,1,\epsilon] = K((\epsilon\gamma +\epsilon + \gamma)/(\gamma +1))$, so 
to be trivial we must have
$\epsilon\gamma +\epsilon + \gamma =\pm 1$. If $\epsilon\gamma +\epsilon + \gamma = 1$,
we get the solutions $\gamma=-2$, $\epsilon=-3$ and $\gamma=-3$, $\epsilon=-2$, 
which correspond  to lines 3-4 of Table 1.

If $\epsilon\gamma +\epsilon + \gamma = -1$, we get the solutions $\epsilon=-1$, 
$\gamma=$arbitrary,
and $\epsilon=$arbitrary, $\gamma=-1$. This gives lines 5-6 of Table 1.

The case $\beta=-1$ is the mirror image of the previous case, so it is not include
in the tables.
The cases when $\gamma$ or $\epsilon$ are $\pm 1$ are similar. 
The case when $\epsilon = 1$
is given in lines 7-10 of Table 1. The case $\gamma=\pm 1 $ is 
symmetric to the case $\beta=\pm 1$.

\begin{center} TABLE 1 \end{center}

\begin{scriptsize}
\begin{center}
\begin{tabular}{|c|c|c|c|c|c|c|}

\hline
 & $\eta$ & $\beta$ & $\alpha$ & $\delta$ & $\epsilon$ & $\gamma$  \\
\hline
 1& 0 & $\pm 1$ & $\pm 1$ & 0 & 0 & $\pm 1$  \\
\hline
 2& 0 &0 & $\pm 1$ & 0 & $\pm 1$ & $\pm 1$ \\
\hline
 & 0 & 1 & $\pm 1$ & 0 & -3 & -2  \\
\hline
 & 0 &1    & $\pm 1$ & 0   & -2  & -3  \\
\hline
& 0 & 1  & $\pm 1$ & 0   & -1  & $\gamma$     \\
\hline
& 0 & 1   & $\pm 1$ & 0   & $\epsilon$  & -1     \\
\hline
& 0 & -2   & $\pm 1$ & 0   &  1 & -3  \\
\hline
& 0 & -3  & $\pm 1$ & 0   &  1 &  -2    \\
\hline
& 0 & $\beta$  & $\pm 1$ & 0   &  1 &  -1    \\
\hline
10 & 0 & -1  &$\pm 1$ & 0   & 1  & $\gamma$     \\
\hline
& 0 & 0 & 1 & 1 & -1 & -2  \\
\hline
 & 0 &0    & 1 & 1   & -2  & -1  \\
\hline
& 0 & 0  & 1 & -1   & $\pm 1-\gamma$  & $\gamma$     \\
\hline
14& 0 & 0   &  -1 &  1   & $\pm 1 -\gamma$  & $\gamma$     \\
\hline

\end{tabular}
\quad
 \begin{tabular}{|c|c|c|c|c|c|c|}
 
 \hline
 & $\eta$ & $\beta$ & $\alpha$ & $\delta$ & $\epsilon$ & $\gamma$  \\
\hline
15& 0 & 0   & -1 & -1   &  1 & 2  \\
\hline
& 0 & 0  & -1 & -1   &  2 &  1    \\
\hline
& 0 & 0  & 1 & -1   &  1 &  -2    \\
\hline
 & 0 & 0  &1 & -2   & 1  & -1     \\
\hline
 & 0 & 0   &  1  & $\pm 1 - \gamma$  & -1  & $\gamma$    \\
\hline
20 & 0 & 0   &  -1  & $\pm 1 - \gamma$  &  1 &  $\gamma$   \\
\hline
& 0 & 0   &  -1  & 1  & -1  &   2  \\
\hline
 & 0 & 0   &  -1  & 2  & -1  &   1  \\
\hline
& 0 & 0   &  1  &  -2 & -1  &  1   \\
\hline
& 0 & 0   &  1  & -1  & -2  &   1  \\
\hline
& 0 & 0   &  1  & $\delta$  & $\pm 1 -\delta$  &  -1   \\
\hline
& 0 & 0   &  -1  &  $\delta$ & $\pm 1 -\delta$  &  1   \\
\hline
& 0 & 0   &  -1  & 2  & 1  &   -1  \\
\hline
28& 0 & 0   &  -1  &  1 & 2  &  -1   \\
\hline

\end{tabular}
\end{center}
\end{scriptsize}

\vfill\eject

Case D. $\beta=0$, $\eta = 0$

In this case the knot is a sum of 2-string tangles. It is made of the Montesinos tangles 
$M(-1/\alpha,-1/\delta)$ and $M(-1/\epsilon,-1/\gamma)$, and well, it is also a Montesinos 
knot. For this to be a trivial knot, one of the two tangle must be a trivial tangle,
and we can assume, because of the symmetry, that the tangle 
$M(-1/\alpha,-1/\delta)$ is trivial. This is trivial only if
$\alpha=\pm 1$ or $\delta=\pm 1$ 

Case D.1. $\alpha=\pm 1$

In this case the knot looks like the Montesinos knot given by 
$M((\mp \delta -1)/\delta, -1/\epsilon,-1/\gamma))$, and to be trivial, one
of the tangles that form it must be an integral tangle.
Suppose first that $R((\mp \delta -1)/-\delta)$ is an integral tangle.
So we have the following three cases:

Case D.1.1. $\alpha=1$, $\delta=1$

In this case the knot is the 2-bridge knot 
$K[\epsilon,2,\gamma]=K((2\epsilon\gamma+\gamma+\epsilon)/(2\epsilon+1))$. 
To be trivial we must have
$2\epsilon\gamma+\gamma+\epsilon=\pm 1$. We get the solutions
$\epsilon=-1$, $\gamma=-2$; $\epsilon=-2$, $\gamma=-1$.
These correspond to lines 11-12 of Table 1.

Case D.1.2. $\alpha=1$, $\delta=-1$; or $\alpha=-1$, $\delta=1$ 

The knot now becomes the 2-bridge knot $K(\epsilon+\gamma)$, which is trivial only if 
$\epsilon+\gamma= \pm 1$. So we get the solutions $\gamma=$arbitrary, 
$\epsilon=\pm 1 -\gamma$. This correspond to lines 13-14 in Table 1.

Case D.1.3. $\alpha=-1$, $\delta=-1$

In this case the knot is the 2-bridge knot 
$K[\epsilon,-2,\gamma]=K((\gamma+\epsilon -2\epsilon\gamma)/(1 -2\epsilon))$. 
To be trivial we must have
$\gamma+\epsilon - 2\epsilon\gamma=\pm 1$. We get the solutions
$\epsilon=1$, $\gamma=2$; $\epsilon=2$, $\gamma=1$.
These correspond to lines 15-16 of Table 1.

The next case inside Case D.1 is to assume that one of the tangles
$R(-1/\epsilon)$ or $R(-1/\gamma)$ is integral. Here the calculation is identical,
and we get lines 17-28 of Table 1. 

Case D.2. $\delta=\pm 1$

This case is symmetric to Case D.1.

Suppose now that just one of the parameters is $0$, say $\eta=0$. In this case the knot looks 
like a sum of two  2-string tangles. See Figure 6.
It is formed by the Montesinos tangles $D(-1/\alpha,-1/\delta)$ and $D(-1/\epsilon,-1/\gamma)$,
which are glued by doing $\beta$ twists. To get the trivial knot, one of the tangles has to be
trivial, and because of the symmetries, we can assume that the tangle $D(-1/\alpha,-1/\delta)$
is trivial. Then $\alpha=\pm 1$, or $\delta=\pm 1$.

\begin{figure}[htbp]
\begin{center}
\includegraphics[width=0.4\linewidth]{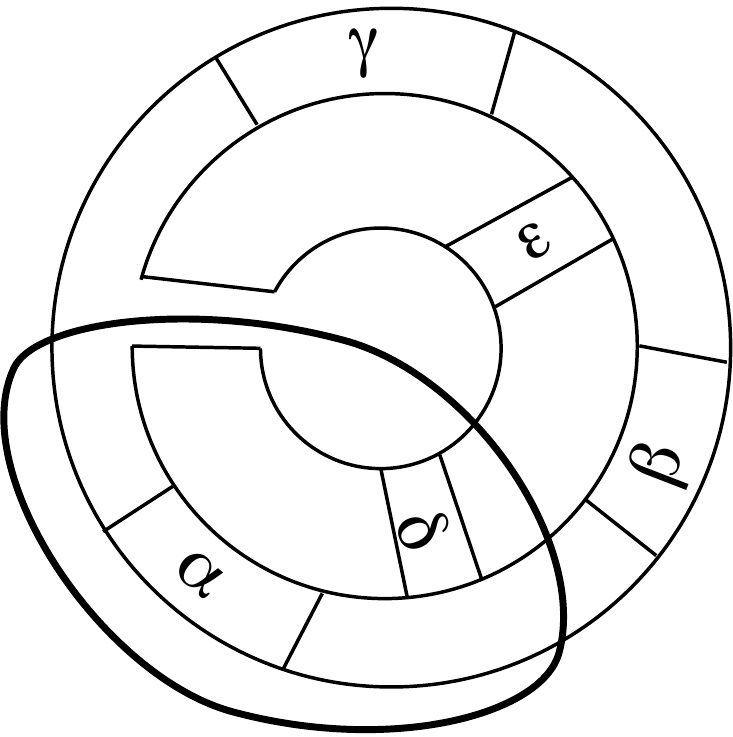}
\caption{}
\label{Figure6}
\end{center}
\end{figure}

Case E. $\alpha=1$

The knot looks like the Montesinos knot 
$M((\delta+1)/(-\beta\delta -\beta-\delta),-1/\epsilon,-1/\gamma)$, and for this to be trivial,
one of the tangles that form it must be integral.

Case E.1. The tangle $R((\delta+1)/(-\beta\delta -\beta-\delta))$ is integral

Then we have $-\beta\delta -\beta-\delta=\pm 1$. We get the solutions:
$\beta=-1$, $\delta$=arbitrary (but to be determined); $\delta=-1$, $\beta=$arbitrary;
$\delta=-2$, $\beta=-3$; $\beta=-2$, $\delta=-3$.

Case E.1.1.  $\beta=-1$, $\delta$=arbitrary (but to be determined)

Now, the knot is the 2-bridge knot $K[\epsilon,-1-\delta,\gamma]=K((-\gamma\epsilon
-\gamma\delta\epsilon +
\gamma +\epsilon)/ (1-\epsilon-\delta\epsilon))$, so for this to be trivial we need that
$-\gamma\epsilon -\gamma\delta\epsilon + \gamma +\epsilon=\pm 1$.
We get the solutions shown in Table 2, lines 1-11.

Case E.1.2. $\delta=-1$, $\beta=$arbitrary

In this case we get the 2-bridge knot $K(\epsilon +\gamma)$, so we get the trivial knot
if $\gamma=$arbitrary, $\epsilon=1-\gamma$; or $\gamma=$arbitrary, $\epsilon=-1-\gamma$.
This gives line 12 in Table 2.

Case E.1.3. $\delta=-2$, $\beta=-3$

In this case the knot is the 2-bridge knot
$K[\epsilon,-1,\gamma]=K((\gamma+\epsilon-\gamma\epsilon)/(1-\epsilon))$, 
so to be trivial we need that
$\gamma+\epsilon-\gamma\epsilon=\pm 1$. We get the solutions
$\gamma=1$, $\epsilon=$arbitrary; $\gamma=$arbitrary, $\epsilon=1$; $\epsilon=2$, $\gamma=3$;
$\epsilon=3$, $\gamma=2$. These correspond to lines 13-16 in Table 2.

Case E.1.4. $\beta=-2$, $\delta=-3$

The knot looks like that 2-bridge knot
$K[\epsilon,-2,\gamma]=K((\gamma+\epsilon-2\gamma\epsilon)/(1-2\epsilon))$, to be trivial we have
$\gamma+\epsilon-2\gamma\epsilon=\pm 1$. We get the solutions
$\epsilon=1$, $\gamma=2$; $\epsilon=2$, $\gamma=1$, which correspond to lines
17-18 in Table 2.

\begin{center} TABLE 2 \end{center}

\begin{scriptsize}
\begin{center}
 \begin{tabular}{|c|c|c|c|c|c|c|}

\hline
 & $\eta$ & $\alpha$ & $\beta$ & $\gamma$ & $\delta$ & $\epsilon$  \\
\hline 
1&0  & 1   &  -1  &$\gamma $  &-1  & $\pm 1-\gamma  $ \\
\hline
 & 0 &  1  & -1   &-2  &-2   &  -3  \\
\hline 
&0  &  1  & -1   &-3  & -2  & -2   \\
\hline 
& 0 &  1  &  -1  &-1  & -3  & -2   \\
\hline 
&0  & 1   & -1   &-2  &-3   &-1    \\
\hline 
&0  & 1   & -1   & -1 & -4  &-1    \\
\hline 
&0  & 1   & -1   &$\gamma $ & -2  & -1   \\
\hline 
& 0 & 1   & -1   & -1 & -2  &$ \epsilon $  \\
\hline 
& 0 & 1   & -1   &1  &1   & 2   \\
\hline 
10&0  &1    & -1   &2  &1   & 1   \\
\hline 
&0  & 1   &-1    &1  & 2  & 1   \\
\hline 
& 0 & 1   &$ \beta $  &$\gamma $ &-1   &$\pm 1-\gamma$    \\
\hline 
& 0 &  1  &   -3 & 1 & -2  &    $\epsilon$ \\
\hline 
& 0 &   1 &  -3  & $\gamma$ &  -2 & 1   \\
\hline 
&  0& 1   &-3    &3  & -2  & 2   \\
\hline 
& 0 & 1   & -3   & 2 &-2   &3    \\
\hline 
&0  & 1   & -2   &2  &-3   &1    \\
\hline 
& 0 & 1   &  -2  & 1 &-3  & 2   \\
\hline 
 & 0 &  1  &  -1  & 1 & 2  & 1   \\
\hline 
20& 0 & 1   &-1    & 2 & 1  & 1   \\
\hline 
 &  0& 1   & -2   &$\gamma$  &$-3-\gamma $  &1    \\
\hline 
& 0 &  1  & -3   & $\gamma $& -2  & 1   \\
\hline 
  & 0   & 1   &$\beta $ &-1   &-2 &1    \\
\hline 
&0  &  1  &$\beta $   &-2  &-1   &1    \\
\hline 
&0  & 1   &-3    &-2  & $\delta $ &1    \\
\hline 
& 0 & 1   &  -2  &$ \gamma$ &$-\gamma -1 $  &  1  \\
\hline 
&0  & 1   &-3    &-3  & -4  &  1  \\
\hline 
&0  & 1   &-3    &-4  & -3  & 1   \\
\hline 
& 0 &  1  &  -4  &-2  &-3   &1    \\
\hline 

30 &0  & 1   &-4    &-3  & -2  &1    \\
\hline 
 &  0&  1  & -5   &-2  &  -2 & 1   \\
\hline 
32& 0 &  1  &$\beta $   &1  &-2   &-1    \\
\hline 
\end{tabular}
\quad
 \begin{tabular}{|c|c|c|c|c|c|c|}
\hline
 & $\eta$ & $\alpha$ & $\beta$ & $\gamma$ & $\delta$ & $\epsilon$  \\
\hline
33& 0 &  1  & -1   &$\gamma $ & -2  & -1   \\
\hline 
& 0 &1    & 1   &3  &2   & -1   \\
\hline 
& 0 & 1   & 2   &2  & 1  & -1   \\
\hline 
& 0 & 1   & 1   &4  & 1  & -1   \\
\hline 
& 0 & 1   & $ \beta$  & 2 & -1  & -1   \\
\hline 
&0  & 1   & 1   & 2 &$\delta $  & -1   \\
\hline 
 &0  &1    & -1   &-1  & -4  & -1   \\
\hline 
40 & 0 &   1 &-2    & -1 &-2   & -1   \\
\hline 
 & 0 &  1  &  -1  &-2  & -3  &-1    \\

\hline 
 & 0 &  1  &  -1  & 1 & 2  & 1   \\
\hline 
& 0 & 1   &-1    & 1 & 1  & 2   \\
\hline 
 &  0& 1   & -2   &1  &$-3-\epsilon $  &$\epsilon$    \\
\hline 
& 0 &  1  & -3   & 1 & -2  & $\epsilon$   \\
\hline 
  & 0   & 1   &$\beta $ &1   &-2 &-1    \\
\hline 
&0  &  1  &$\beta $   & 1 &-1   &-2    \\
\hline 
&0  & 1   &-3    &1  & $\delta $ &-2    \\
\hline 
& 0 & 1   &  -2  &1 &$-\epsilon -1 $  &  $\epsilon$  \\
\hline 
50&0  & 1   &-3    &1  & -4  &  -3  \\
\hline 
&0  & 1   &-3    &1  & -3  & -4   \\
\hline 
& 0 &  1  &  -4  &1  &-3   &-2    \\
\hline 

 &0  & 1   &-4    & 1 & -2  &-3    \\
\hline 
 &  0&  1  & -5   &1  &  -2 & -2   \\
\hline 
& 0 &  1  &$\beta $   &-1  &-2   &1    \\
\hline 
& 0 &  1  & -1   &-1 & -2  & $\epsilon$   \\
\hline 
& 0 &1    & 1   &-1  &2   & 3   \\
\hline 
& 0 & 1   & 2   &-1  & 1  & 2   \\
\hline 
& 0 & 1   & 1   &-1  & 1  & 4   \\
\hline 
60& 0 & 1   & $ \beta$  & -1 & -1  & 2   \\
\hline 
&0  & 1   & 1   & -1 &$\delta $  & 2   \\
\hline 
 &0  &1    & -1   &-1  & -4  & -1   \\
\hline 
 & 0 &   1 &-2    & -1 &-2   & -1   \\
\hline 
64 & 0 &  1  &  -1  &-1  & -3  &-2    \\
\hline

\end{tabular}
\end{center}
\end{scriptsize}

Case E.2. The tangle $R(-1/\epsilon)$ is integral

Case E.2.1. $\epsilon=1$

The knot is the 2-bridge knot $K[\delta,1,\beta,1,\gamma]=K((\gamma\beta\delta
+\gamma\beta+2\delta\gamma+\beta\delta +
\gamma +\beta +\delta)/(\beta\delta + \beta +2\delta +1))$. To be trivial, the numerator
must be $\pm 1$. We get the solutions shown in lines 19-31 of Table 2.

Case E.2.2. $\epsilon=-1$

The knot is the 2-bridge knot $K[\delta,1,\beta,-1,\gamma]=K((-\gamma\beta\delta - \gamma\beta
+\gamma+\beta\delta +\beta+\delta)/(-\beta\delta+\beta+1))$. Again to be trivial, the numerator must
be $\pm 1$. We get the solutions shown in lines 32-41 of Table 2.

Case E.3. The tangle $R(-1/\gamma)$ is integral

The analysis is identical to the case E.2, just interchanging $\epsilon$ and
$\gamma$. We get the solutions shown in lines 42-64 of Table 2.

Case F. $\alpha=-1$

This is the mirror image of Case E), so it is not shown in the tables.

Case G. $\delta=\pm 1$

In this case $\alpha$ and $\delta$ can be interchanged by a reflection on the hexatangle,
giving solutions equivalent to the previously found.
\qed

\section{When a parameter is $1$ and the other is $-1$}
\label{sec:paremeter1}

{\bf Theorem 1.4}\ {\it Suppose that $\alpha$, $\beta$, $\gamma$,
$\epsilon$ and $\eta$ are not $0$, $\delta=-1$ and $\eta=1$.
$\mathcal {H} (\alpha,\beta,\gamma,-1, \epsilon,1)$ is the trivial
knot if and only the parameters are as in Table 3, up to symmetries.}

\textbf{Proof}
In this case the knot looks like a Montesinos knot, see Figure 7.
In fact, it is the Montesinos knot 
$M(\alpha/(1-\epsilon\alpha), \gamma/(\gamma +1),\beta/(1-\beta))$.
This knot can be trivial only if one of the rational tangles that form it is 
an integral tangle, and this happens only if the denominators of
the fractions are $\pm 1$. So if $1-\epsilon\alpha=\pm 1$, we have the
cases: $\alpha=1$, $\epsilon=2$; $\alpha=-1$, $\epsilon=-2$; 
$\alpha=2$, $\epsilon=1$; $\alpha=-2$, $\epsilon=-1$. If $\gamma+1=\pm 1$,
then $\gamma=-2$. If $1-\beta=\pm 1$, then $\beta=2$ (remember that we
are assuming that none of the parameters is $0$).

\begin{figure}[htbp]
\begin{center}
\includegraphics[width=0.8\linewidth]{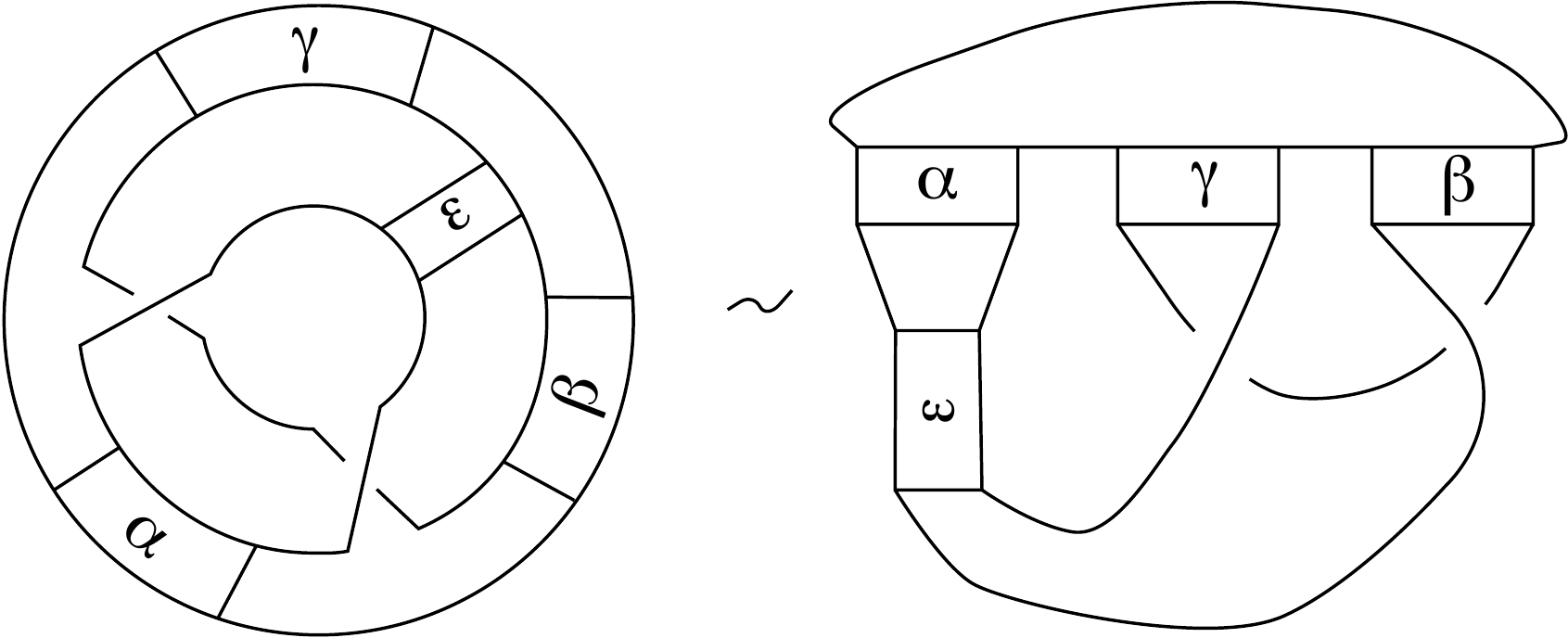}
\caption{}
\label{Figure7}
\end{center}
\end{figure}

We have the following cases:

Case A. $\alpha=1$, $\epsilon=2$

In this case the knot can be isotoped so that it looks like the numerator
of a rational tangle, and by computing the continued fraction, we see that it 
is the 2-bridge knot $K[-\beta,\gamma+2]=K((1-\beta\gamma-2\beta)/\beta)$.
For this to be trivial, it is needed that the knot is of the form $1/n$,
that is,  $1-\beta\gamma-2\beta=\pm 1$. We get the solutions:
$\gamma=-2$, $\beta=$arbitrary; $\beta=-1$, $\gamma=-4$;
$\beta=2$, $\gamma=-1$; $\beta=-2$, $\gamma=-3$. These correspond to lines
1-4 of Table 3.

Case B. $\alpha=-1$, $\epsilon=-2$

In this case the knot is the 2-bridge knot
$K[2-\beta,\gamma]=K((2\gamma-\beta\gamma+1)/(2-\beta))$. This is trivial if and only if
$2\gamma-\beta\gamma+1=\pm 1$. 
We get the solutions:
$\beta=2$, $\gamma=$arbitrary; $\beta=4$, $\gamma=1$;
$\beta=3$, $\gamma=2$; $\beta=1$, $\gamma=-2$. These correspond to lines
5-8 of Table 3.

\begin{center} TABLE 3 \end{center}

\begin{scriptsize}
\begin{center}
 \begin{tabular}{|c|c|c|c|c|c|c|}

\hline
 & $\delta$ & $\eta$ & $\alpha$ & $\beta$ & $\gamma$ & $\epsilon$  \\
\hline
1& -1 & 1   &  1  &$ \beta$ &-2   & 2   \\
\hline 
&-1  & 1   &1    & -1 &-4   &2    \\
\hline 
& -1 & 1   & 1   &2  & -1  & 2   \\
\hline 
&-1  & 1   & 1   & -2 &-3   & 2   \\
\hline 
& -1 & 1   &-1    &2  &$ \gamma$  & -2   \\
\hline 
& -1 &   1 &-1    & 4 &1   & -2   \\
\hline 
&-1  & 1   & -1   &1  & -2  &  -2  \\
\hline 
& -1 & 1   &  -1  &  3&2   &-2    \\
\hline 
&-1  & 1   & 2   & 2 & -1  &1    \\
\hline 
10 & -1 &  1  &2    & -1 & -2  & 1   \\
\hline 
&-1  & 1   & 2   & 1 &-2   &1    \\
\hline 
&-1  &1    &-2    &1  &-2   &  -1  \\
\hline 
& -1 &1    &-2    &2  & 1  & -1   \\
\hline 
& -1 & 1   & -2   & 2 & -1  &-1    \\
\hline 
& -1 & 1   &  $\alpha $ &$-\alpha +1 $ &-2   &1    \\
\hline 
&-1  &1    &  2  & 3 &-2   &  3  \\
\hline 
&-1  & 1   & 2   & 5 & -2  & 2	   \\
\hline 
&-1  & 1   &3    & 4 & -2  & 2   \\
\hline 
& -1 & 1   &  1  & 3 &-2   & 4   \\
\hline 
20 & -1 &  1  & 1   &4  &-2   & 3   \\
\hline

\end{tabular}
\quad
\begin{tabular}{|c|c|c|c|c|c|c|}

\hline
 & $\delta$ & $\eta$ & $\alpha$ & $\beta$ & $\gamma$ & $\epsilon$  \\
\hline 
21& -1 &  1  &  -1  &2  &-2   & $\epsilon $  \\
\hline 
&-1  & 1   & 1   &$\beta$  & -2  & 2   \\
\hline
&-1 & 1   &$\alpha$   &3  &-2   & 2   \\
\hline
& -1 & 1   & $\alpha $  &$3-\alpha $ &-2   & 1   \\
\hline
& -1 & 1   &-1   & 1 &-2   & -2   \\
\hline
&-1  & 1   &  -2  &1  & -2  &  -1  \\
\hline
& -1 & 1   &  1  & 2 &-2   & 3   \\
\hline
&-1  & 1   &  1  & 2 &  -2 & 4   \\
\hline
 &-1  & 1   & 1   &2  & -2  &$ \epsilon $  \\
\hline
30& -1 & 1   &$ \alpha$   &2  &$-1-\alpha $  & -1   \\
\hline
&-1  &1    & -2   &2  &  -3 & -3   \\
\hline
&-1  & 1   & -2   &2  &-5   &-2    \\
\hline
& -1 &  1  & -3   &2 & -4  & -2   \\
\hline
&-1  & 1   &   -1 &2  &-3   & -4   \\
\hline
&-1  & 1   &  -1  &2  &-4   & -3   \\
\hline
&-1  &1    &$\alpha $   &1  & -3  &-2    \\
\hline
&-1  & 1   &   -1 &2  &$\gamma $  &-2    \\
\hline
& -1 &  1  & $ \alpha $ &2  &$-3-\alpha$   &-1    \\
\hline
&-1  & 1   & 1   &2  & -1  & 2   \\
\hline
40& -1 & 1   & 2   & 2 &-1   &  1  \\
\hline

 \end{tabular}
\end{center}
\end{scriptsize}

Case C. $\alpha=2$, $\epsilon=1$

In this case the knot is the 2-bridge knot
$K[\gamma+2,-2,\beta]=K((2\gamma\beta + 3\beta -\gamma-2)/(2\gamma +3))$. 
For this to be trivial
we need that $2\gamma\beta + 3\beta -\gamma-2=\pm 1$, and this is possible only
in the following cases: $\gamma=-1$, $\beta=2$; $\gamma=-2$, $\beta=-1$;
$\gamma=-2$, $\beta=1$. These correspond to lines 9-11 of Table 3.

Case D. $\alpha=-2$, $\epsilon=-1$

In this case the knot is the 2-bridge knot
$K[\gamma,2,\beta-2]=K((2\gamma\beta + \beta -3\gamma-2)/(2\gamma +1))$. 
For this to be trivial
we need that $2\gamma\beta + \beta -3\gamma-2=\pm 1$, and this is possible only
in the following cases: $\gamma=-2$, $\beta=1$; $\gamma=$, $\beta=2$;
$\gamma=-1$, $\beta=2$. These correspond to lines 12-14 of Table 3.

Case E. $\gamma=-2$

In this case the knot is the 2-bridge knot
$K[\alpha,-\epsilon,2,-1,\beta]=$
$K((\alpha\beta\epsilon-\alpha\beta-2\alpha\epsilon + \alpha -\beta+2)/(\alpha\epsilon-\alpha-1))$.
For this to be trivial we need that
$\alpha\beta\epsilon-\alpha\beta-2\alpha\epsilon + \alpha -\beta+2=\pm 1$.
A careful calculation shows that this is possible only for the cases shown in
lines 15-29 of Table 3.

Case F. $\beta=2$

In this case the knot is the 2-bridge knot 
$K[\alpha,-\epsilon,-2,1,\gamma]=$
$K((\epsilon\alpha\gamma + \alpha\gamma - \gamma + 2\epsilon\alpha-2+\alpha)/
(\epsilon\alpha+\alpha-1))$. Again, this is trivial if and only if 
$\epsilon\alpha\gamma + \alpha\gamma - \gamma + 2\epsilon\alpha-2+\alpha=\pm 1$.
A careful calculation show that this is possible only for the cases
shown in lines 21 y 29-40 of Table 3.
\qed

\section{The reduction Lemmas}
\label{sec:mainlemmas}

In this section we prove the following theorem.

{\bf Theorem 1.3} \ {\it Suppose that all of 
$\alpha$, $\beta$ $\gamma$, $\delta$,
$\epsilon$ and $\eta$ are different from $0$. If 
$\mathcal {H} (\alpha,\beta,\gamma,\delta, \epsilon,\eta)$ is the trivial
knot, then there is a pair of adjacent boxes, 
say $\delta$ and $\eta$,
so that $\delta=-1$ and $\eta=1$.}

This is proved in several steps.

\subsection{The first reduction}
\label{subsec:firstred}

Let $V_1$, $V_2$, $V_3$, $V_4$ be solid tori. Let $M_1= V_1 \cup_A V_2$,
where $V_1$ and $V_2$ are glued along an annulus 
$A\subset \partial V_i$, $i=1,2$, and suppose that $A$ goes at least
twice longitudinally on each solid tori. Let $M_2= V_3 \cup_B V_4$,
where $V_3$ and $V_4$ are glued along an annulus 
$B\subset \partial V_i$, $i=3,4$, and suppose that $B$ goes at least
twice longitudinally on each solid tori. So $M_1$ and $M_2$ are of the
form $D(a,b)$, i.e., a Seifert fibered spaces over the disk with two cone points.
Let $\lambda_i$ be a fiber on $\partial M_i$ of its Seifert fibering.
In the special case that the annulus $A$ or the annulus $B$ goes exactly twice
longitudinally on each solid tori, the manifold $M_i$ is also a Seifert fibered space
over the M\"obius band without cone points; let $\lambda_i'$ be a fiber on $\partial M_i$
of this fibering.
 Let $N=M_1\cup M_2$, glued along their boundary,
and suppose that $\lambda_1$, and $\lambda_1'$, if this is the case, are not identified 
to  curves isotopic to $\lambda_2$ or $\lambda_2'$.
Then any of the corresponding fibers have geometric intersection number $\geq 1$
in the torus $T=\partial M_1=\partial M_2$.
So $N$ is a graph manifold, containing the incompressible torus $T$
which divide it into two Seifert fibered spaces, but $N$ is not a
Seifert fibered space. Let $k_1$ ($k_2$) be an embedded arc in the
annulus $A$ (resp. $B$) joining points on different components of
$\partial A$ (resp. $\partial B$). Assume that 
$\partial k_1 = \partial k_2$, and let $k=k_1\cup k_2$. So $k$ is a
knot in $N$. Finally let $M=N -int\, \eta(k)$. So $M$ is a
compact 3-manifold with a torus boundary component.

\begin {LEM} \label{LE:hyp}The manifold $M$ is irreducible, atoroidal, and not a
Seifert  fibered space, hence it is hyperbolic. $M$ has a Dehn filling
which produces the toroidal manifold $N$.
\end{LEM}

\textbf{Proof} Let $T$ be the essential torus in $N$, i.e.,
$T=\partial M_1 = \partial M_2$, and let 
$\hat T = T-int\,\eta (k)$. So $\hat T$ is a twice punctured torus
properly embedded in $M$. We show first that $\hat T$ is incompressible in
$M$. Suppose $D$ is a compression disk for $\hat T$, then, say, $D$ is
contained in $M_1-int\, \eta (k_1)$. But as $T$ is incompressible in
$M_1$, $\partial D$ is inessential in $T$, so $\partial D$ must bound
a disk on $T$ which contains the points $k_1 \cap T$, but this would
imply that the arc $k_1$ is contained in a 3-ball, which is not
possible. 

$M$ is irreducible, for it is the union of two irreducible manifolds
glued along an incompressible surface. 

Suppose $R$ is an incompressible torus in $M$. As $\hat T$ and $R$
are incompressible, they can be isotoped so that their intersection consists
of simple closed curves which are essential in both surfaces.
These divide $R$ into a collection of annuli $A_1,\dots,A_n$.
Suppose $A_1$ lies on $M_1-int\, \eta(k_1)$. Look now at the intersections
between $A_1$ and $E_1=A-int\, \eta(k_1)$. Note that $E_1$ is a disk, so trivial 
curves of intersection can be eliminated. If there is an arc of intersection
whose endpoints are in the same component of $\partial A$ and also in the same
component of $\partial A_1$, then this can be eliminated by an
isotopy of $A_1$. If there is an arc of intersection whose endpoints are in the 
same component of $\partial A$ but in different components of $\partial A_1$, 
then $A_1$ is $\partial$-compressible and it follows that it is parallel onto $\hat T$. 
If this happens, then by pushing $A_1$ onto $M_2$, the number of curves of 
intersection between $\hat T$ and $R$ is reduced. Also, there are no arcs
of intersection whose endpoints are in the same component of $\partial A_1$ but 
different components of $A$, for $A$ would be $\partial$-compressible. 
So any arc of intersection has endpoints in different components of $\partial A$, 
and in different components of $\partial A_1$. These  arcs of intersection cut
$A_1$ into squares, and it follows that there are 3 possibilities for $A_1$,  
either it is  parallel to $\partial \eta(k_1)$, or it is parallel to $A$, or $A$ goes
twice longitudinally on $V_1$ and $V_2$ and $A_1$ is formed by the union of two squares,
one lying in $V_1$ and the other in $V_2$. Note that in the last two cases, $\partial A_1$
is isotopic to a fiber in a Seifert fibration of $M_1$.
Now look at $A_2$, which lies in $M_2 -int\, \eta(k_2)$. A similar analysis can be done
in this case. If $A_2$ is parallel to $\hat T$, again an isotopy reduces intersections. 
By construction, $\partial A_2$ cannot be isotopic to a fiber in a Seifert fibration of $M_2$.
So the only possibility left is that $A_1$ is parallel to $\partial \eta (k_1)$, and
that $A_2$ is parallel to $\partial \eta (k_2)$. This implies that
$R$ is peripheral, that is, isotopic to $\partial M$.

As $M$ is atoroidal, if it is a Seifert fibered space, it must be of
the form $D(a,b)$, but then a Dehn filling of it cannot be a graph
manifold.
\qed

Note that if we remove from $M$ the core of one of the solid tori in $M_1$ and the core of one of the solid tori in $M_2$, then the resulting manifold $M'$ has 3 tori as boundary components, and the same proof shows that it is hyperbolic.

As an approximation to Theorem 1.3, we first show the following.

\begin{LEM} \label{LE:6C}If $\mathcal {H} (\alpha, \beta,\gamma,\delta, \epsilon,\eta)$
is the trivial knot, where none of $\alpha$, $\beta$, $\gamma$, $\delta$, $\epsilon$,
$\eta$ is $0$, then one of the following cases must occur:

\item 1. There is a pair of opposite boxes, say $\eta$ and $\beta$, so that all the other 
boxes are different from $\pm 1$. Furthermore there are the following cases, up to 
symmetries and mirror images:

a)\quad  $\beta =1$, $\eta =2$

b)\quad  $\beta =-1$, $\eta =2$

c)\quad  $\beta =1$, $\eta =1$ 

d)\quad $\beta =-1$, $\eta =1$

\item 2. There is a pair of adjacent boxes, say $\delta$ and $\eta$,
which are 1 or -1. So we have the following cases, up to  mirror images:

e)\quad $\delta =1$, $\eta =1$

f)\quad $\delta =-1$, $\eta =1$

\end{LEM}

\textbf{Proof} Let $\tilde {\mathcal {H}} (\alpha, \beta,\gamma,\delta, \epsilon,\eta)$ denote
the double branched cover of \break
$\mathcal {H} (\alpha,\beta,\gamma,\delta, \epsilon,\eta)$.
Suppose that there is a pair of opposite boxes, say $\eta$ and $\beta$, 
so that the other boxes $\vert \alpha \vert, \vert \gamma
\vert, \vert \delta \vert, \vert \epsilon \vert$ are $\geq 2$. Note
that in $\mathcal {H} (\alpha, 0,\gamma,\delta, \epsilon,\eta)$ there is a
sphere $S$ which decomposes it as a sum of prime tangles, something similar
to Figure 6. A lift of this sphere
in $\tilde {\mathcal {H}} (\alpha, 0,\gamma,\delta, \epsilon,\eta)$ is an
incompressible torus. 

Note that 
$\tilde {\mathcal {H}} (\alpha, \beta,\gamma,\delta, \epsilon,\eta)$ can be
identified with $\mathcal {L}(1/\eta,1/\beta,1/\epsilon,-\gamma,-\alpha,-\delta)$
(making  a $2\pi/3$-rotation on Figure 4 to match Figure 3).
So $\tilde {\mathcal {H}} (\alpha, 0,\gamma,\delta, \epsilon,\eta)=
\mathcal {L}(1/\eta,1/0,1/\epsilon,-\gamma,-\alpha,-\delta)$, which is depicted
in Figure 8. From it, it is easy to see that there is a torus $T$, which
divides the surgered manifold into two Seifert fibered spaces, one is given by
a solid torus containing knots with framings $-\gamma$ and $-\alpha$,
and the other by a solid torus containing knots with framings $-\delta$ and $1/\epsilon$.
These are glued in a twisted way, given by the knot with framing $1/\eta$. This gluing
ensures that the Seifert fibers of one side are not identified to the fibers on the other side,
except possibly if all of $\alpha$, $\gamma$, $\delta$, $\epsilon$, are $\pm 2$, but this case
is not relevant in our case, for if it happens then the hexatangle produces a link of 2 or more components, not a trivial knot.
It follows that $\tilde {\mathcal {H}} (\alpha, *,\gamma,\delta, \epsilon,\eta)=
\mathcal {L}(1/\eta,*,1/\epsilon,-\gamma,-\alpha,-\delta)$  is a manifold
as $M$ in Lemma \ref{LE:hyp}. The knot to be 
removed from $\mathcal {L}(1/\eta,1/0,1/\epsilon,-\gamma,-\alpha,-\delta)$ 
to get $M$ is shown with dotted lines in Figure 8. 
Now it follows from Lemma \ref{LE:hyp} that 
$\tilde {\mathcal {H}} (\alpha, *,\gamma,\delta, \epsilon,\eta)$
is a hyperbolic manifold.

\begin{figure}[htbp]
\begin{center}
\includegraphics[width=0.4\linewidth]{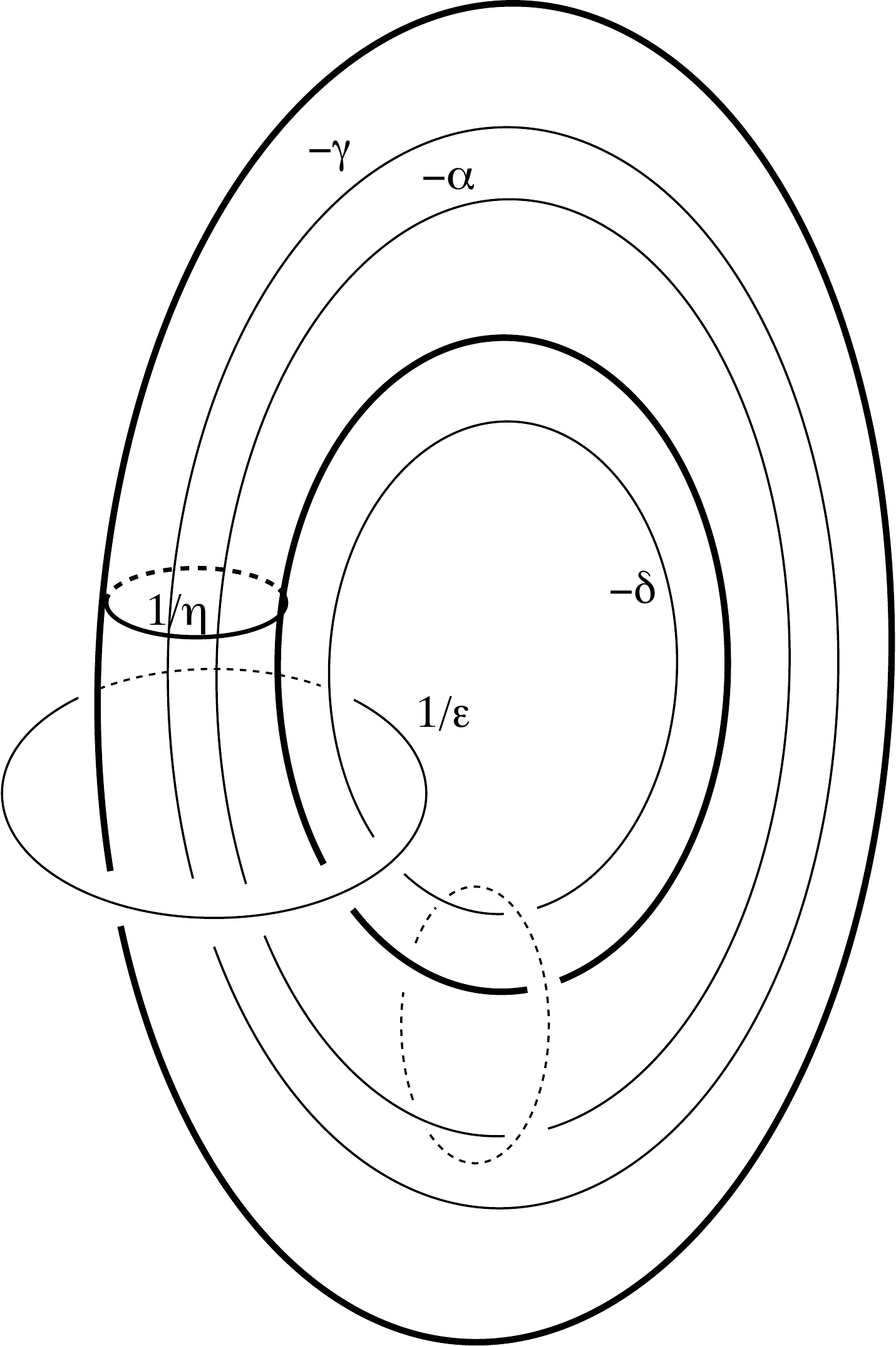}
\caption{}
\label{Figure8}
\end{center}
\end{figure}

So if $\tilde {\mathcal {H}} (\alpha, \beta,\gamma,\delta, \epsilon,\eta)$ is the 3-sphere, 
we must have that
$\vert \beta\vert \leq 2$ (by \cite{GordonL2}, or \cite{Eudave1}), so $\beta = \pm 1$ or $\beta= \pm 2$.
In any case, $\tilde {\mathcal {H}} (\alpha, \beta,\gamma,\delta,
\epsilon,*)$ is hyperbolic, again by Lemma 1.1. So by the same argument
we get that $\eta$ is $\pm 1$ or $\pm 2$. So we get Case (1) of the Lemma,
except if $\vert \beta\vert =\vert \eta\vert =2$. If this happens, then take another
pair of opposite boxes, say $\alpha$ and $\epsilon$,  and repeat the
argument. Again we get Case (1), or $\vert \alpha\vert =\vert \epsilon\vert =2$.
But note that ${\mathcal {H}} (\pm 2, \pm 2,\gamma,\delta, \pm 2,\pm 2)$
is a link of two or more components, so we must have Case (1) of the Lemma

Now, if for any pair of opposite boxes one of the remaining parameters is $\pm 1$,
then there is a pair of adjacent boxes, both of which are $\pm 1$.
So we have Case (2) of the Lemma.
\qed

These six cases are shown in Figure 9.

\begin{figure}[htbp]
\begin{center}
\includegraphics[width=0.9\linewidth]{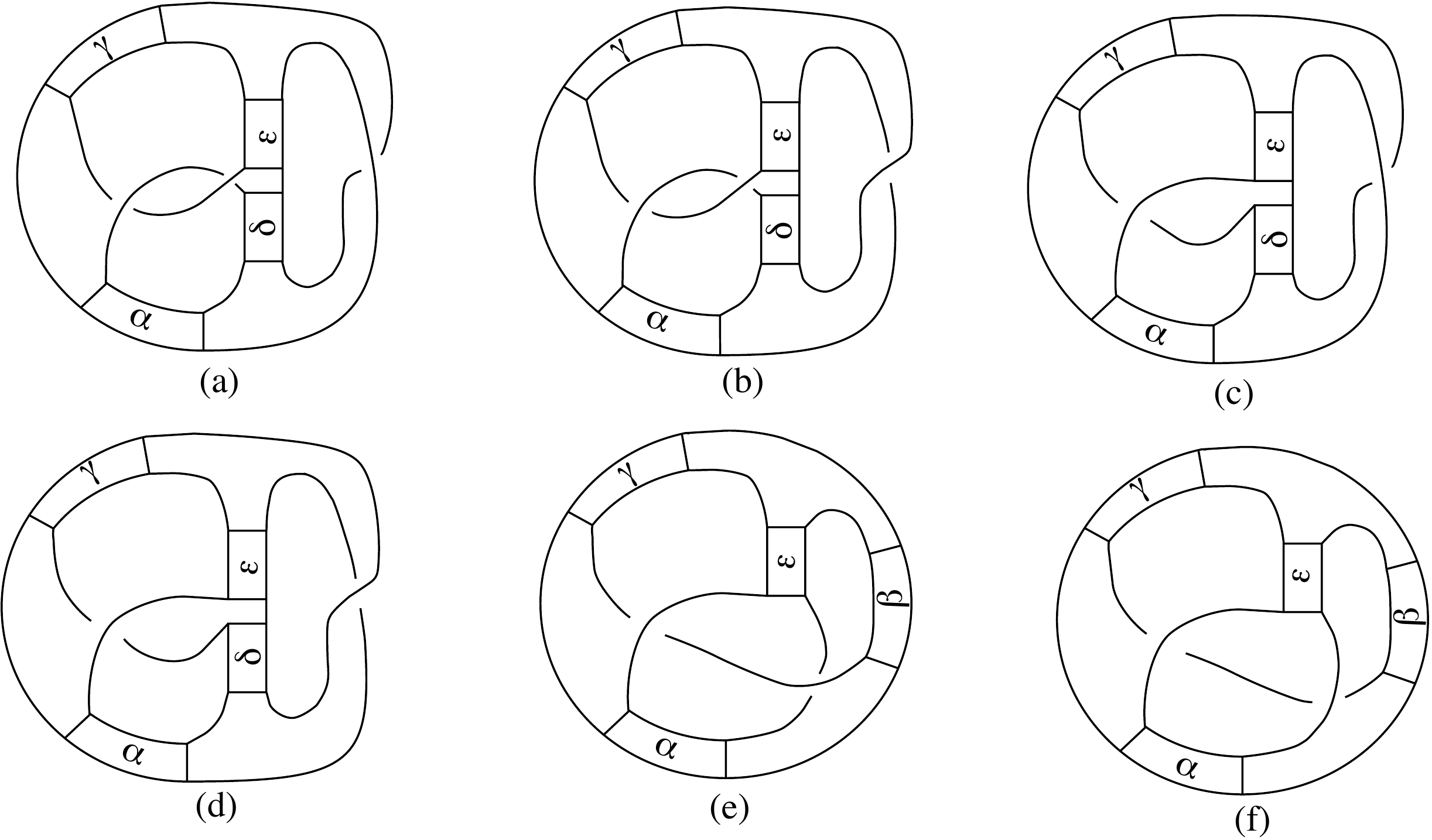}
\caption{}
\label{Figure9}
\end{center}
\end{figure}

\subsection{The knots $\mathcal {Q}_1 (\alpha,\gamma,\delta,\epsilon)$}

Let $\mathcal {Q}_1 (\alpha,\gamma,\delta,\epsilon)= 
\mathcal {H} (\alpha,1,\gamma,\delta,\epsilon,2)$,
as shown in Figure 9(a).

In this section we prove the following:

\begin{PROP}\label{PR:Q1} If all of $\vert \alpha \vert$, $\vert \gamma \vert$, 
$\vert \delta \vert$, $\vert \epsilon \vert$ are
$\geq 2$, then $\mathcal {Q}_1 (\alpha,\gamma,\delta,\epsilon)$ cannot be the trivial knot.
\end{PROP}

\textbf{Proof} Suppose that $\mathcal {Q}_1 (\alpha,\gamma,\delta,\epsilon)$
is the trivial knot.
By Lemma 4.1.1, \break $\tilde {\mathcal {H}}(\alpha,1,\gamma,\delta,\epsilon,*)$ is a hyperbolic manifold,
in fact the exterior of a hyperbolic knot in $S^3$ by hypothesis. It also has a toroidal filling,
corresponding to $\eta=0$, which is at distance $2$ from the filling $\eta=2$ that produces the
3-sphere. So, by the main result of \cite{GordonL3}, we have that 
$\tilde {\mathcal {H}}(\alpha,1,\gamma,\delta,\epsilon,*)$ is the exterior of one of the knots
$k(\ell,m,n,p)$ constructed in \cite{Eudave2}. It follows that $\tilde {\mathcal {H}}(\alpha,1,\gamma,\delta,\epsilon,0)$ double branch covers a EM-knot $K(\ell,m,n,p)$ as defined in \cite{GordonL4}, but also it double branch covers $\mathcal {H}(\alpha,1,\gamma,\delta,\epsilon,0)$, so Theorem 3.4 of \cite{GordonL4} implies that $\mathcal {H}(\alpha,1,\gamma,\delta,\epsilon,0)$ must be one of the EM-knots. 
Let $S(*,*,*,*)$ be the tangle defined in \cite{GordonL4}, Figures 2.3, 2.4.
Note that $\mathcal {H}(\alpha,1,\gamma,\delta,\epsilon,0)=
S(-1/\alpha,(-\delta-1)/\delta,-1/\epsilon,(-\gamma-1)/\gamma)$.

It follows from \cite{Eudave3}, 5.4, or \cite{GordonL4}, 3.1 that the EM-knot
$K(\ell,m,n,p)$ is the same as $S(\alpha',\beta',\gamma',\delta')$ , where $\alpha'$, $\beta'$, $\gamma'$,
$\delta'$ are as follows:

$$p=0 : \ \ \alpha'= -\frac{1}{\ell},\  \beta'=\frac{m}{\ell m-1},\ \gamma'=\frac{2mn+1-m-n}{4mn-2m+1},\
\delta'=-1/2$$

$$n=0 : \ \alpha'=-\frac{1}{\ell},\  \beta'=\frac{2mp-m-p}{\ell(2mp-m-p)-2p+1},
\ \gamma'=\frac{m-1}{2m-1},\ \delta'=-1/2$$

By Lemma 2.2 of \cite{GordonL4}, and by the symmetries of 
$S(-1/\alpha,(-\delta-1)/\delta,-1/\epsilon,\break (-\gamma-1)/\gamma)$, we can assume that the Montesinos tangles $M(\alpha',\beta')$ and \break 
$M(-1/\alpha,(-\delta-1)/\delta)$ are equivalent as marked tangles, and also are the tangles
$M(\gamma',\delta')$ and $M(-1/\epsilon,(\gamma-1)/\gamma)$.

Case A. $p=0$

The numerator of  $M(-1/\ell, m/(\ell m-1))$ is the trivial knot, and the numerator of 
$M(-1/\alpha,(-\delta-1)/\delta)$ is the 2-bridge knot
$K((\alpha\delta +\alpha +\delta)/(\alpha+1))$, which can be trivial only if $\alpha=-2$, $\delta=-3$, or
$\alpha=-3$, $\delta=-2$ (or $\alpha=\pm 1, 0$). So $M(-1/\ell,m/(\ell m -1))$
is the tangle $M(1/2,-1/3)=M(-1/2,2/3)$ or the tangle $M(-1/3,1/2)$. This is possible only if $\ell=-3$, $m=-1$, or $\ell=-2,m=-2$, or $\ell=2, m=2$.

Case A.1. $m=-1$

Looking at the other pair of tangles in the decomposition, we get that $M(-1/\epsilon, \break (-\gamma-1)/\gamma)=M((-\epsilon-1)/\epsilon,-1/\gamma)$ is equivalent to $M((3n-2)/(4n-3),-1/2)$. There are two cases.

Case A.1.1. $\gamma=2$ and $(-\epsilon-1)/\epsilon=(3n-2)/(4n-3)$

It follows that $(3n-2)+(4n-3)=\pm 1$, which is impossible for an integral $n$.

Case A.1.2. $\gamma=-2$ and $-1/\epsilon=(3n-2)/(4n-3)$

In this case we have $3n-2=\pm 1$, which is possible only if $n=1$, but then $\epsilon=-1$, which is not possible by hypothesis.

Case A.2. $m=-2$

Looking again at the other pair of tangles in the decomposition, we get that $M(-1/\epsilon,(-\gamma-1)/\gamma)$ is equivalent to $M((3-5n)/(5-8n),-1/2)$. There are two cases as above, and a similar argument show that they are not possible.

Case A.3. $m=2$

We get that $M(-1/\epsilon,(-\gamma-1)/\gamma)$ is equivalent to $M((3n-1)/(8n-3),-1/2)$.
An argument as before shows that this is not possible.

Case B. $n=0$

The tangle $M((m-1)/(2m-1),-1/2)$ is equivalent to the tangle $M(-1/\epsilon,\break (-\gamma-1)/\gamma)=M((-\epsilon-1)\epsilon,-1/\gamma)$. There are two cases:

Case B.1. $\gamma=2$ and $(-\epsilon-1)/\epsilon = (m-1)/(2m-1)$

We have that $(m-1)+(2m-1)=\pm 1$, which is possible only if $m=1$. However in this case
we have that $\epsilon=1$.

Case B.2. $\gamma=-2$ and $-1/\epsilon = (m-1)/(2m-1)$

So $m-1=\pm 1$, which implies that $m=2$ and $\epsilon=-3$.

Look at the other pair of tangles. Having $m=2$, we get that $M(-1/\alpha, (-\delta-1)/\delta)$ is 
equivalent to $M(-1/\ell,(3p-2)/(3p\ell-2\ell-2p+1))$. We have the following cases:

Case B.2.1. $(-\delta-1)/\delta = (3p-2)/(3p\ell-2\ell-2p+1)$

Then $p+3p\ell -2\ell -1=\pm 1$. A calculation shows that this is possible only if
$\ell=\pm 1, 0$, but then $\alpha=\pm 1, 0$.

Case B.2.2. $\ell =2$, $\alpha=-2$, and $M(1/2,(-\delta-1)/\delta)=M(-1/2,-1/\delta)= M(-1/2, (3p-2)/(4p-3))$

Then $-1/\delta = (3p-2)/(4p-3)$. So $3p-2=\pm 1$, which is possible only if $p=1$, but in this case
$\delta=-1$.

Case B.2. $\ell =-2$, $\alpha=2$, and $M(-1/2,(-\delta-1)/\delta)=M(1/2,(-2\delta-1)/\delta)= M(1/2, (3p-2)/(-8p+5))$

Then $(-2\delta-1)/\delta= (3p-2)/(-8p+5))$, which implies that  $(3p-2)+2(-8p+5)=\pm 1$, which is clearly not possible. \qed

\subsection{The knots  
$\mathcal {Q}_2 (\alpha,\gamma,\delta,\epsilon)$}

Let $\mathcal {Q}_2 (\alpha,\gamma,\delta,\epsilon)= 
\mathcal {H} (\alpha,-1,\gamma,\delta,\epsilon,2)$,
as shown in Figure 9(b).

In this section we prove the following:

\begin{PROP} \label{PR:Q2}If all of $\vert \alpha \vert$, $\vert \gamma \vert$, 
$\vert \delta \vert$,
$\vert \epsilon \vert$ are
$\geq 2$, then $\mathcal {Q}_2 (\alpha,\gamma,\delta,\epsilon)$ cannot be the trivial knot.
\end{PROP}

\textbf{Proof} The proof follows the same lines as Proposition \ref{PR:Q1}, just
note that  
$\mathcal {H}(\alpha,-1,\gamma,\delta,\epsilon,*)$ 
$=S(-1/\alpha,(\delta-1)/\delta,-1/\epsilon,(\gamma-1)/\gamma)$.
Assume that as marked tangles $M(-1/\alpha,(\delta-1)/\delta)=M(\alpha',\beta')$, and
$M(-1/\epsilon,(\gamma-1)/\gamma)=M(\gamma',\delta')$.

Case A. $p=0$

Then $M(-1/\ell, m/(m\ell-1))= M(-1/\alpha,(\delta-1)/\delta)$. The numerator of the first tangle is the trivial knot, and the numerator of the second tangle is the 2-bridge knot 
$K((\alpha+\delta-\alpha\delta)/(1-\alpha))$, which is the trivial knot only if $\alpha=2$, $\delta=2$ or $\alpha=3$, $\delta=2$. So $M(-1/\alpha,(\delta-1)/\delta)$ is the tangle $M(-1/2,2/3)=M(1/2,-1/3)$ or the tangle $M(-1/3,1/2)$. This is possible only if $\ell=2$, $m=2$, or $\ell=-2$, $m=1$, or
$\ell=3$, $m=1$.

Case A.1. $m=2$

Looking at the other pair of tangle we get that
$M((3n-1)/(8n-3),-1/2)= M(-1/\epsilon,(\gamma-1)/\gamma) = M((\epsilon-1)/\epsilon,-1/\gamma)
= M((2\epsilon-1)/\epsilon,(-1-\gamma)/\gamma)$.
We have the cases:

Case A.1.1. $\gamma=2$, and then $(3n-1)/(8n-3)=(\epsilon-1)/\epsilon$

Case A.1.2. $\gamma=-2$, and then $(3n-1)/(8n-3)=(2\epsilon-1)/\epsilon$

A calculation shows that these cases are not possible.

Case A.2. $m=1$

We have that $M(n/(4n-1),-1/2)= M(-1/\epsilon,(\gamma-1)/\gamma) = M((\epsilon-1)/\epsilon,-1/\gamma) =M((2\epsilon-1)/\epsilon,(-1-\gamma)/\gamma)$. We have the cases:

Case A.1.1. $\gamma=2$, and then $n/(4n-1)=(\epsilon-1)/\epsilon$

Case A.1.2. $\gamma=-2$, and then $n/(4n-1)=(2\epsilon-1)/\epsilon$

A calculation shows that these cases are not possible.

Case B. $n=0$

In this case we have $M((m-1)/(2m-1),-1/2) = M(-1/\epsilon,(\gamma-1)/\gamma) = 
 M((\epsilon-1)/\epsilon,-1/\gamma) = M((2\epsilon-1)/\epsilon, (-\gamma-1)/\gamma)$.

Case B.1. $\gamma=-2$, and then $(m-1)/(2m-1)=(2\epsilon-1)/\epsilon$

A simple calculation shows that this case is not possible.

Case B.2. $\gamma=2$, and then $(m-1)/(2m-1)=(\epsilon-1)/\epsilon$

This is possible only if $m=\pm 1$, but if $m=1$ then $\epsilon=1$.

Case B.2.1. $m=-1$

Looking at the other pair of tangles we get that $M(-1/\ell,(-3p+1)/(-3p\ell + \ell -2p+1))= M(-1/\alpha,(\delta-1)/\delta)$. We have the following cases:

Case B.2.1. $(\delta-1)/\delta = (-3p+1)/(-3p\ell+\ell-2p+1)$

Then $-p+3p\ell -\ell =\pm 1$. A calculation shows that this is possible only if
$\ell=\pm 1, 0$, but then $\alpha=\pm 1, 0$.

Case B.2.2. $\ell =2$, $\alpha=-2$, and $M(1/2,(\delta-1)/\delta)=M(-1/2,(2\delta-1)/\delta)= M(-1/2, (3p-1)/(8p-3))$

Then $(2\delta-1)/\delta = (3p-1)/(8p-3)$. So $(3p-1)-2(8p-3)=\pm 1$, which is clearly not possible.

Case B.2.3. $\ell =-2$, $\alpha=2$, and $M(-1/2,(\delta-1)/\delta)=M(1/2,-1/\delta)= M(1/2, (-3p+1)/(4p-1))$

Then $-1/\delta= (-3p+1)/(4p-1))$, which implies that  $-3p+1=\pm 1$, which is clearly not possible. 
\qed

\subsection{The knots $\mathcal {Q}_3 (\alpha,\gamma,\delta,\epsilon) $}

Let $\mathcal {Q}_3 (\alpha,\gamma,\delta,\epsilon)= 
\mathcal {H} (\alpha,1,\gamma,\delta,\epsilon,1)$, as in Figure 9(c), and let  
$\tilde{\mathcal {Q}}_3 (\alpha,\gamma,\delta,\epsilon)$ be its double branched cover.
In this section we prove the following:

\begin{PROP} \label{PR:Q3} If all of $\vert \alpha \vert$, $\vert \gamma \vert$, 
$\vert \delta \vert$,
$\vert \epsilon \vert$ are
$\geq 2$, then $\mathcal {Q}_3 (\alpha,\gamma,\delta,\epsilon)$ 
cannot be the trivial knot.
\end{PROP}

\textbf{Proof} 
Suppose that $\mathcal {Q}_3(\alpha,\gamma,\delta,\epsilon)$ is the trivial knot.

\begin{CLAIM} Either $\alpha=-2$ or $\epsilon =-2$.
\end{CLAIM}

\textbf{Proof} Consider the tangle $\mathcal {Q}_3(\alpha,\gamma,\delta,*)$, then
$\tilde {\mathcal {Q}}_3(\alpha,\gamma,\delta,*)$ is the exterior of a knot in $S^3$.
Note that $\mathcal {Q}_3(\alpha,\gamma,\delta,-1)$ looks like a composite knot.
In fact, it is the connected sum of two-bridge knots
$K(\alpha +1)\# K((\delta\gamma-1)/ \gamma)$,
which will be in fact composite, unless $\alpha =-2$ (or one of $\delta$, 
$\gamma$ is $0$ or $\pm 1$).
Suppose then that $\alpha\not= -2$. As the knot is composite, its double 
branched cover is reducible, and then the corresponding surgery must be at
distance 1 from $\epsilon$ \cite{GordonL1}, so $\epsilon=-2$.
\qed

Suppose then that $\epsilon=-2$. By symmetry, the other case is identical.
Consider then $\mathcal {Q}_3(\alpha,\gamma,\delta,-2)$.
Note that this looks like the Montesinos knot
$M(-1/(\gamma +1),-1/ (\delta+1),\alpha/(1+2\alpha))$, and for this to be trivial,
one of the rational tangles that form it must be an integral tangle,
which is possible only if $\gamma=-2$ or $\delta=-2$ (or $\alpha=-1$, which
is not considered by our hypothesis).

Suppose first that $\gamma=-2$. Then the knot $\mathcal {Q}_3(\alpha,-2,\delta,-2)$
is the two bridge knot $K[-\alpha,-2,-1,\delta+1]= K((3\alpha\delta + \delta +\alpha)/(3\alpha+1))$.
To be trivial, we need that $3\alpha\delta + \delta + \alpha = \pm 1$. We see
by inspection that this is not possible, unless one of $\delta$, $\alpha$ were 
$0$ or $\pm 1$.

If now we suppose that $\delta =-2$, the same arguments produce a contradiction.
\qed

\subsection{ The knots $\mathcal {Q}_4 (\alpha,\gamma,\delta,\epsilon)$}

Let $\mathcal {Q}_4 (\alpha,\gamma,\delta,\epsilon)= 
\mathcal {H} (\alpha,1,\gamma,\delta,\epsilon,-1)$, as in Figure 9(d),
and let  $\tilde {\mathcal {Q}}_4
(\alpha,\gamma,\delta,\epsilon)$ its double branched cover.
In this section we prove the following:

\begin{PROP} \label{PR:Q4} If all of $\vert \alpha \vert$, $\vert \gamma \vert$, 
$\vert \delta \vert$, $\vert \epsilon \vert$ are
$\geq 2$, then $\mathcal {Q}_4 (\alpha,\gamma,\delta,\epsilon)$ 
cannot be the trivial knot.
\end{PROP}

\textbf{Proof} Suppose that $\mathcal {Q}_4(\alpha,\gamma,\delta,\epsilon)$ is the trivial knot.
Then $\tilde{\mathcal {Q}}_4(\alpha,\gamma,\delta,*)$ is the exterior of a knot in $S^3$.

\begin{CLAIM} The tangle $\mathcal {Q}_4(\alpha,\gamma,\delta,*)$
is trivial.
\end{CLAIM}

\textbf{Proof} Note that
$\mathcal {Q}_4(\alpha,\gamma,\delta,0)$ is a 2-bridge knot. 
Then the knot $\tilde {\mathcal {Q}}_4(\alpha,\gamma,\delta,*)$ has a Dehn surgery
producing a lens space. If the knot is not a torus knot, nor a trivial knot,
then such surgery must be a distance 1 from $\epsilon$ [CGLS], which is not
possible in our case, for $\vert \epsilon \vert \geq 2$. 
So the knot must be a torus knot or the trivial knot.
If it is a torus knot, then it must have a reducible surgery at distance
one from the lens space surgery and at distance one from $\epsilon$, 
so the possibilities are $1/0$, $1$ or $-1$. Note that these fillings produces
knots that look like Montesinos knots, in fact 
$\mathcal {Q}_4(\alpha,\gamma,\delta,1/0)= M(-1/\alpha,-1/(1+\delta),1/(1-\gamma))$,
so this can be composite only if
$\alpha=0$, $\gamma=1$ or $\delta =-1$. We have that
$\mathcal {Q}_4(\alpha,\gamma,\delta,1)= M(\delta/(1-\delta\gamma),-1/2,1/(1-\alpha))$,
which can be composite only if $\alpha=1$, or $\delta\gamma=1$.
$\mathcal {Q}_4(\alpha,\gamma,\delta,-1)=M(\gamma/(1-\delta\gamma),1/2,-1/(1+\alpha))$ 
which can be composite only if $\delta\gamma =1$, of $\alpha=-1$.
All these possibilities are not possible in our case,
so the tangle $\mathcal {Q}_4(\alpha,\gamma,\delta,*)$
is trivial.  
\qed

As the tangle ${\mathcal {Q}}_4(\alpha,\gamma,\delta,*)$ is trivial, any filling
of it must give a 2-bridge knot. Note that 
$\mathcal {Q}_4(\alpha,\gamma,\delta,-1)$ looks like the Montesinos knot
$M(\gamma/(1-\delta\gamma)$, $1/2,-1/(\alpha+1))$, and for this to be a 
2-bridge knot, one of the rational tangle that form it must be an 
integral tangle, and this is possible only if $\alpha=-2$ 
(or one of $\gamma$, $\delta$ is $0$ or $\pm 1$).

Suppose then that $\alpha=-2$. Note that $\mathcal {Q}_4(-2,\gamma,\delta,1)$
looks like the Montesinos knot $M(\delta/(1-\delta\gamma),-1/2,1/3)$, and for this to be a 2-bridge 
knot, we need that $\delta\gamma=2$, so one of $\delta$ or $\gamma$ must be $\pm 1$,
which is not possible.
\qed

\subsection{The knots $\mathcal {Q}_5 (\alpha,\beta,\gamma,\epsilon)$ }

Let $\mathcal {Q}_5 (\alpha,\beta, \gamma,\epsilon)= 
\mathcal {H} (\alpha,\beta,\gamma,1,\epsilon,1)$, as in Figure 9(e),
and let  $\tilde {\mathcal {Q}}_5
(\alpha,\beta,\gamma,\epsilon)$ its double branched cover.
In this section we prove the following:

\begin{PROP} \label{PR:Q5}
 If all of $\vert \alpha \vert$, $\vert \beta \vert$, 
$\vert \gamma \vert$, $\vert \epsilon \vert$ are
$\geq 1$, and $\alpha \not= -1$, $\beta \not= -1$, $\gamma \not= -1 $ 
and $\epsilon \not= -1$ then 
$\mathcal {Q}_5 (\alpha,\beta,\gamma,\epsilon)$ cannot be the trivial knot.
\end{PROP}

\begin {LEM} \label{LE:Q5}
Suppose that $\alpha$, $\beta$ and $\gamma$
are as in Proposition  \ref{PR:Q5}, and that  $\epsilon =1$. Then 
$\mathcal {Q}_5 (\alpha,\beta,\gamma,1)$ is not the trivial knot.
$\mathcal {Q}_5 (\alpha,\beta,\gamma,1)$ is a prime link, except if, up to symmetries,
$\alpha=-2$ and $\beta=-\gamma$.
\end{LEM}

\textbf{Proof} Suppose first that $\alpha=-2$. Note that $\mathcal {Q}_5 (-2,\beta,\gamma,1)$
is the Montesinos link $M(-1/(\beta +\gamma),-2/3,1/2)$. This can be composite
only if $\beta+\gamma=0$, and in fact in this case we get the connected sum of a
trefoil knot with the Hopf link. The knot can be trivial only if $\beta+\gamma=\pm 1$.
If $\beta+\gamma= 1$, we get the two-bridge knot $K(7/2)$, and if
$\beta+\gamma=- 1$, we get the 2-bridge knot $K(-5)$, so no trivial knot is obtained.

Suppose then that $\alpha$, $\beta$, $\gamma$ are $\not= -2$.
Note that the knot $\mathcal {Q}_5 (\alpha,\beta,\gamma,1)$ is a closed 3-braid
around an axis perpendicular to the plane which passes through the
triangle formed by $\delta$, $\eta$, $\epsilon$. In fact, it is the closed braid \break
$\sigma_1^{-\alpha}\sigma_2\sigma_1^{-\beta}\sigma_2\sigma_1^{-\gamma}\sigma_2$.
We will apply the classification of links which are closed 3-braids \cite{BirmanM}, to conclude that
our knot is not trivial nor composite. To do that we first find the Schreier unique
representative of the conjugacy class of the braid (see \cite{BirmanM}, \S 7). We have the following cases:

Case A. $\alpha$, $\beta$ and $\gamma$ are positive

In this case the braid $\sigma_1^{-\alpha}\sigma_2\sigma_1^{-\beta}\sigma_2\sigma_1^{-\gamma}\sigma_2$
is already the Schreier unique representative of its conjugacy class.

\vfill\eject

Case B. $\alpha$ is negative, and $\beta$, $\gamma$ are positive

In this case, we calculate the Schreier unique
representative of its conjugacy class to be 
$C\sigma_1^{-\gamma-1}\sigma_2^{-\alpha-2}\sigma_1^{-\beta-1}\sigma_2$, where
$C=(\sigma_1\sigma_2\sigma_1)^2 = (\sigma_1\sigma_2)^3$.

Case C. $\alpha$, $\beta$ are negative, and $\gamma$ is positive

In this case we get $C^2\sigma_1^{-\gamma-2}\sigma_2^{-\alpha-3}\sigma_1^{-1}\sigma_2^{-\beta-3}$.

Case D. $\alpha$, $\beta$ and $\gamma$ are negative

In this case we get 
$C^3\sigma_1^{-1}\sigma_2^{-\alpha-4}\sigma_1^{-1}\sigma_2^{-\beta-4}\sigma_1^{-1}\sigma_2^{-\gamma-4}$,
in the generic case. There are several special cases. If $\alpha=-3$, $\beta=-3$, 
we get $C^2\sigma_1^{-\gamma-3}$. If $\alpha=-3$, $\beta=-4$, $\gamma=-4$,
we get $C^2\sigma_1\sigma_2$. If $\alpha=-3$, $\beta=-4$, $\gamma=-5$,
we get $C^2\sigma_1\sigma_2\sigma_1$. If $\alpha=-3$, $\beta=-4$, $\gamma=-6$,
we get $C^2\sigma_1\sigma_2\sigma_1\sigma_2$. If $\alpha=-3$, $\beta=-4$, $\gamma<-6$,
we get $C^3\sigma_1^{-1}\sigma_2^{-\gamma-7}$. If $\alpha=-3$, $\beta<-4$,
we get $C^3\sigma_1^{-1}\sigma_2^{-\beta-5}\sigma_1^{-1}\sigma_2^{-\gamma-5}$.

The remaining cases are identical to the given ones, because of the symmetries of the hexatangle.
The trivial knot has three conjugacy classes of 3-braid representatives, namely: 
$\sigma_1\sigma_2$, $\sigma_1^{-1}\sigma_2$ and $C^{-1}(\sigma_1\sigma_2)^2$. None of these braids was
obtained in our calculation, so we conclude that $\mathcal {Q}_5 (\alpha,\beta,\gamma,1)$ is not
the trivial knot.

A composite link which is a closed 3-braid has as Schreier representative of the conjugacy class
of a 3-braid representing it, one of the following: $\sigma_1^{-u}\sigma_2^v$, where $u\geq v\geq 2$, or
$C^{-1}\sigma_1^{-u}\sigma_2\sigma_1^{-v}\sigma_2$, where $u\geq v\geq 0$ (see \cite{BirmanM},\S 7). None of these
braids was obtained in our calculation, so we conclude that $\mathcal {Q}_5 (\alpha,\beta,\gamma,1)$
is not a composite link, except as said before, if one of $\alpha$, $\beta$ or $\gamma$ is $=-2$.
\qed

\textbf{Proof of Proposition \ref{PR:Q5}}
Note that $\mathcal {Q}_5 (\alpha,\beta, \gamma,0)$ is a 2-bridge knot. By \cite{CGLS},
either $\epsilon=\pm 1$, or $\tilde {\mathcal {Q}}_5 (\alpha,\beta, \gamma,*)$ is 
the exterior of a torus or a trivial knot.  If $\epsilon=1$, we are done by Lemma \ref{LE:Q5}.
So suppose first that it is the
exterior of a torus knot. Then there is a filling of 
$\tilde {\mathcal {Q}}_5 (\alpha,\beta, \gamma,*)$, say $\epsilon_1\in Q\cup\{1/0\}$, that produces 
a reducible manifold, and then $\mathcal {Q}_5 (\alpha,\beta, \gamma,\epsilon_1)$
will be a composite link. The slope $\epsilon_1$ must be a distance 1 from 
both, the slope $0$ and $\epsilon$. There are the following possibilities:
A) $\epsilon_1=1/0$; B) $\epsilon_1=1$, $\epsilon=2$; 
C) $\epsilon_1=-1$, $\epsilon=-2$; D) $\epsilon_1= 1/2$, $\epsilon=1$; 
E) $\epsilon_1=-1/2$, $\epsilon=-1$. If Case D happens, we finish, and
Case E does not happen by hypothesis. So we have Cases A, B or C.

Case A. $\epsilon_1=1/0$

The link $\mathcal {Q}_5 (\alpha,\beta, \gamma,1/0)$ looks like the Montesinos knot
$M(-1/\alpha,-1/2,$ \break   $ -1/(\beta+\gamma))$,
but it must be composite, then one of the rational tangle that form it must be
the tangle $1/0$, so we have that $\alpha=0$, which is not possible, or $\beta=-\gamma$. 
So assume that $\beta=-\gamma$. The reducible manifold will be 
$L(2,1)\# L(\alpha,1)$, but remember that $pq$-Dehn surgery on the torus knot
$T_{p,q}$ produces the lens space $L(p,q)\# L(q,p)$, so we must have that $p=2$, and
that $L(q,2)$ is homeomorphic to $L(\alpha,1)$, but this is possible only
if $\alpha=\pm 3$.

Suppose then that $\alpha=\pm 3$. Any surgery on $\tilde {\mathcal {Q}}_5 (\alpha,\beta, \gamma,*)$ at distance one 
from $1/0$ must produce a lens space.
So $\tilde {\mathcal {Q}}_5 (\pm 3 ,\beta, -\beta,-1)$ must be a lens space,
and then $\mathcal {Q}_5 (\pm 3, \beta, -\beta,-1)$ is a 2-bridge link. 
But note that it looks like the Montesinos link $M(3/4,-\beta/(\beta+1),-\beta/(\beta-1))$
or $M(3/2,-\beta/(\beta +1),-\beta/(\beta-1))$, depending if $\alpha=3$ or $\alpha=-3$,
 and then to be a
2-bridge knot we must have $\beta=2$ or $-2$, but by symmetry these cases are identical,
so suppose $\beta=2$. Doing the filling $\epsilon=1$, we also have to get
a lens space, so  $\mathcal {Q}_5 (3,2, -2,1)$ must be a 2-bridge knot, 
but it is the Montesinos knot $M(1/2,-2/3,-1/5)$, which is not a 2-bridge knot,
so $\alpha\not=3$.

Note that $\tilde {\mathcal {Q}}_5(-3,2,-2,*)$ is in fact the exterior of the trefoil knot,
but  $\mathcal {Q}_5(-3,2,-2,0)$ is the trivial knot, so $\mathcal {Q}_5(-3,2,-2,\epsilon)$ cannot
be trivial for $\epsilon\not= 0$.

Case B. $\epsilon_1=1$, $\epsilon=2$

So the knot $\mathcal {Q}_5 (\alpha,\beta, \gamma,1)$ must be composite, but this is not
possible by Lemma \ref{LE:Q5}, unless one of $\alpha$, $\beta$ or $\gamma$ is $-2$.
Suppose first that $\alpha=-2$, and $\gamma=-\beta$. Note that
$\mathcal {Q}_5 (-2,\beta, -\beta,2)$ is not trivial, it is the connected sum of
the figure eight knot and the Hopf link.

Suppose now that $\gamma=-2$, and that $\beta=-\alpha$.
Any filling at distance 1 from $\epsilon_1$ must produce a lens space.
$\mathcal {Q}_5 (\alpha,-\alpha, -2,1/0)$ looks like the Montesinos link
$M(-1/\alpha,-1/2,1/(\alpha+2))$, so for this to be a 2-bridge knot we must
have $\alpha=\pm 1, -3$. The case $\alpha=-1$ is not considered by hypothesis. 
Suppose $\alpha=1$. The knot $\mathcal {Q}_5 (1,-1, -2,2)$ is not trivial, in fact,
it is the 2-bridge knot $K(13/5)$. Suppose now that $\alpha=-3$.
The knot $\mathcal {Q}_5 (-3,3, -2,2)$ is not trivial,
it is the 2-bridge knot $K(13/5)$. The case $\beta=\pm 2$, $\gamma=-\alpha$ is symmetric
to the previous case.

Case C. $\epsilon_1=-1$, $\epsilon=-2$

So the knot $\mathcal {Q}_5 (\alpha,\beta, \gamma,-1)$ must be composite. Note that it 
looks like the Montesinos knot $M(\gamma/(1-\gamma),\beta/(1-\beta),\alpha/(1+\alpha))$,
 so to be composite one of the tangles that 
form it must be $1/0$, so we have that either $\beta=1$ or $\gamma=1$ (or $\alpha=-1$,
but in this case we finish). Note also that both cases are symmetric,
so we can assume that $\beta=1$.

The knot $\mathcal {Q}_5 (\alpha,1, \gamma,1/0)$ must be a 2-bridge knot, but it looks
like the Montesinos knot $M(-1/\alpha,-1/2,-1/(1+\gamma))$, so to be a 2-bridge knot
we must have $\alpha=\pm 1$ or $\gamma=-2$.
If $\alpha=-1$ we finish, so we have the other two cases.

\vfill\eject

Case C.1. $\alpha=1$

The knot $\mathcal {Q}_5 (1,1, \gamma,-2)$ looks like the Montesinos knot 
$M(1/3,-1/2,-1/(\gamma+1))$, it is a 2-bridge knot only if $\gamma=0$ or $\gamma=-2$,
but if $\gamma=-2$ then it is the 2-bridge knot $K(5)$, which is not trivial.

Case C.2. $\gamma=-2$

Then $\mathcal {Q}_5 (\alpha,1,-2,-2)$ is the 2-bridge knot $K((-4\alpha-1)/\alpha))$,
which cannot be trivial if $\alpha$ is integral.

So we have shown that $\tilde {\mathcal {Q}}_5 (\alpha,\beta, \gamma,*)$ cannot be
the exterior of a torus knot. Suppose now that it is the exterior of
the trivial knot. Then all the links obtained from 
$\mathcal {Q}_5 (\alpha,\beta, \gamma,*)$ must be 2-bridge links. 
$\mathcal {Q}_5 (\alpha,\beta, \gamma,-1)$ looks like the Montesinos knot
$M(\gamma/(1-\gamma),\beta/(1-\beta),\alpha/(1+\alpha))$, so to be a 2-bridge link
we must have $\beta=2$, $\gamma=2$, or $\alpha=-2$.

Case F.1. $\beta=2$

$\mathcal {Q}_5 (\alpha,2, \gamma,1/0)$ looks like the Montesinos knot
$M(-1/\alpha,-1/2,-1/(2+\gamma))$, but it  must be a 2-bridge knot,
so we have $\alpha=\pm 1$, $\gamma=-1$ or $\gamma=-3$. 
If $\alpha=-1$ or $\gamma=-1$ we finish.

Case F.1.1. $\alpha=1$

$\mathcal {Q}_5 (1,2, \gamma,-1/2)$ looks like the sum of the two Montesinos tangles 
$T(-1/\gamma,$ $-1/2)$ and $T(-1/3,1/2)$, so to be a 2-bridge knot we need that $\gamma=\pm 1$.
If $\gamma=-1$ we finish. If $\gamma=1$, then $\mathcal {Q}_5 (1,2,1,-1/2)$ is the Montesinos
knot $M(2/3,-1/3,1/2)$, which is not a 2-bridge knot.

Case F.1.2. $\gamma=-3$

$\mathcal {Q}_5 (\alpha,2, -3,-1/2)$ looks like the sum of 
two Montesinos tangles, $T(-1/2,1/3)$ and $T((\alpha+1)/(\alpha+2),-1/2)$, 
so to be a 2-bridge knot  we need that $\alpha=- 1$ or
$\alpha=-3$. But $\mathcal {Q}_5 (-3,2, -3,-1/2)$ is the Montesinos knot 
$M(-2/3,-1/2,1/3)$, which is not a 2-bridge knot.

Case F.2. $\gamma=2$. This is symmetric to the case F.1.

Case F.3. $\alpha=-2$

$\mathcal {Q}_5 (-2,\beta, \gamma,1/0)$ looks like the Montesinos knot
$M(1/2,-1/2,-1/(\beta+\gamma))$, which is a 2-bridge knot only when $\beta+\gamma=\pm 1$.
If $\beta+\gamma=1$, then $\mathcal {Q}_5 (-2,\beta, \gamma,\epsilon)$ is the 2-bridge knot
$K(-(4\epsilon+3)/4)$, which is trivial only if $\epsilon=-1$. If $\beta+\gamma=1$, then 
$\mathcal {Q}_5 (-2,\beta, \gamma,\epsilon)$ is the 2-bridge knot
$K((4\epsilon+1)/4)$, which is trivial only if $\epsilon=0$.  
\qed

\textbf{Proof of Theorem \ref{THM:quasi-main}}
Suppose that $\mathcal {H} (\alpha,\beta,\gamma,\delta, \epsilon,\eta)$ is the trivial knot and that
all of $\alpha$, $\beta$ $\gamma$, $\delta$, $\epsilon$ and $\eta$ are different from $0$.
There are two cases which exclude each other: (A) there is a pair of opposite boxes so that all of the other boxes have two or more crossings, or (B) there is a pair of adjacent boxes each with a single crossing. In case A, Lemma \ref{LE:6C} shows that there are 4 possible cases, up to symmetries, these are: (a) $\beta =1$, $\eta =2$; (b) $\beta =-1$, $\eta =2$; (c) $\beta =1$, $\eta =1$; 
(d) $\beta =-1$, $\eta =1$. Propositions \ref{PR:Q1}, \ref{PR:Q2}, \ref{PR:Q3} and \ref{PR:Q4} show that none of these cases produce the trivial knot.
For case B, there are two possibilities up to symmetries, either $\delta=1$ and $\eta=1$ and none of the other parameters is $-1$, or
$\delta=-1$ and $\eta=1$. The first case cannot produce the trivial knot by Proposition \ref{PR:Q5},
so we must have that $\delta=-1$ and $\eta=1$.
\qed

\section {The pure 3-braids with trivial surgeries} 

In this section we return to the original problem of determining which small closed pure 3-braids
produce $S^3$ by surgery. All information is contained in Tables 1, 2, 3. Any entry in the tables produces many braids with a trivial surgery, because of the symmetries. We will not reproduce all such braids here. Now we will look at some specific and interesting examples of 3-braids with a trivial surgery.

Let $\mathcal {L}$ be the link shown in Figure 2. As before, we indicate surgeries on this link
by $\mathcal {L}(1/e_1,1/f_1,1/e,m,n,p)$, as indicated in Figure 2, which implicitly is
giving an order to the components of the link. Note that a surgery $\mathcal {L}(1/0,*,*,*,*,*)$,
or a surgery $\mathcal {L}(1,*,-1,*,*,*)$ produces a non-hyperbolic manifold. So, any entry in the tables
corresponds to surgery on some non-hyperbolic braid. However, we do get many hyperbolic braids
from such tables.

\begin {PROP} There exists infinitely many hyperbolic small 3-braids which have
a non-trivial surgery producing the 3-sphere.
\end{PROP}

\textbf{Proof} Take the solution of line 23 in Table 3. We have $\alpha$, $\beta=3$, $\gamma=-2$,
$\delta=-1$, $\eta=1$, $\epsilon=2$ (i.e, $\alpha$ takes any value). 
This is equivalent to the solution $\alpha=3$, $\beta=-2$,
$\gamma$, $\delta=2$, $\eta=-1$, $\epsilon=1$. The double branched cover of
$\mathcal {H}(3,-2,\gamma,2,1,-1)$ is $\mathcal {L}(1/2,1/\gamma,-1,-3,2,-1)$, obtained by identifying Figures 3 and 4 without any transformation. Consider the closed pure 3-braid $\mathcal {L}(1/2,1/\gamma,-1,*,*,*)$. The proof of Lemma \ref{LE:hyp} can be applied to show that this braid is hyperbolic if $\gamma\not= \pm 1,\ 0$, by the remark made just after the proof of the Lemma. We apply the proof so that the fifth component of $\mathcal {L}$, i.e., the one which covers the box $B$ is the one we are filling to get a toroidal manifold.
So, the closed braids $\sigma_1^4\sigma_2^{2\gamma}(\sigma_2\sigma_1\sigma_2)^{-2}$ are all hyperbolic if $\vert\gamma\vert\geq 2$. Each of these braids have a non-trivial surgery producing 
$S^3$.
\qed 

We have shown that the closed braid $\hat \beta =\widehat{\sigma_1^4\sigma_2^{2\gamma}(\sigma_2\sigma_1\sigma_2)^{-2}}$ is hyperbolic if $\vert\gamma\vert\geq 2$, and according to SnapPea \cite{Weeks}, it is also hyperbolic for $\gamma=-1$. Each of these braids have a non-trivial surgery producing 
$S^3$; by adjusting the surgery coefficients, we have that the surgery that produces $S^3$ is 
$(-4,1-\gamma,-\gamma)$. 

\textbf{Acknowledgement} We are grateful to Francisco Gonz\'alez-Acu\~na and to the referee for valuable suggestions. This research was partially supported by PAPIIT--UNAM grant IN115105.

\end{document}